\documentclass[11pt]{article}

\usepackage{authblk}

\RequirePackage[OT1]{fontenc}
\RequirePackage{amsthm,amsmath}
\RequirePackage[numbers]{natbib}
\RequirePackage[colorlinks,citecolor=blue,urlcolor=blue]{hyperref}

\usepackage{latexsym, epsfig, amssymb, amsmath, amsthm, graphicx, mathrsfs}
\usepackage[usenames,dvipsnames]{color}
\usepackage{anysize}

\marginsize{1.2in}{1in}{1.35in}{1in} %

\theoremstyle{plain}
\newtheorem{theorem}{Theorem}
\newtheorem{definition}{Definition}
\newtheorem{lemma}{Lemma}
\newtheorem{assumption}{Condition}

\definecolor{gray}{rgb}{0.7,0.7,0.7}

\newcommand{\ba}{\mbox{\bf a}}
\newcommand{\vertiii}[1]{{\left\vert\kern-0.25ex\left\vert\kern-0.25ex\left\vert #1
\right\vert\kern-0.25ex\right\vert\kern-0.25ex\right\vert}}

\newcommand{\bb}{\mbox{\bf b}}

\newcommand{\br}{\mbox{\bf r}}

\newcommand{\bu}{\mbox{\bf u}}
\newcommand{\bv}{\mbox{\bf v}}
\newcommand{\bw}{\mbox{\bf w}}
\newcommand{\bx}{\mbox{\bf x}}
\newcommand{\by}{\mbox{\bf y}}
\newcommand{\bz}{\mbox{\bf z}}
\newcommand{\bA}{\mbox{\bf A}}
\newcommand{\bB}{\mbox{\bf B}}
\newcommand{\bC}{\mbox{\bf C}}
\newcommand{\bD}{\mbox{\bf D}}

\newcommand{\bG}{\mbox{\bf G}}
\newcommand{\bM}{\mbox{\bf M}}
\newcommand{\bN}{\mbox{\bf N}}

\newcommand{\mR}{\mathbb{R}}

\newcommand{\bH}{\mbox{\bf H}}
\newcommand{\bI}{\mbox{\bf I}}

\newcommand{\bR}{\mbox{\bf R}}
\newcommand{\bq}{\mbox{\bf q}}
\newcommand{\bS}{\mbox{\bf S}}

\newcommand{\bV}{\mbox{\bf V}}
\newcommand{\bW}{\mbox{\bf W}}
\newcommand{\bX}{\mbox{\bf X}}

\newcommand{\bY}{\mbox{\bf Y}}
\newcommand{\bZ}{\mbox{\bf Z}}

\newcommand{\bone}{\mbox{\bf 1}}
\newcommand{\bzero}{\mbox{\bf 0}}
\newcommand{\bveps}{\mbox{\boldmath $\varepsilon$}}

\newcommand{\bbeta}{\mbox{\boldmath $\beta$}}
\newcommand{\bdelta}{\mbox{\boldmath $\delta$}}
\newcommand{\btheta}{\mbox{\boldmath $\theta$}}

\newcommand{\bPsi}{\mbox{\boldmath $\Psi$}}

\newcommand{\bmu}{\mbox{\boldmath $\mu$}}

\newcommand{\bSigma}{\mbox{\boldmath $\Sigma$}}

\newcommand{\hB}{\widehat \bB}
\newcommand{\hA}{\widehat \bA}
\newcommand{\hH}{\widehat \bH}

\newcommand{\hbbeta}{\widehat\bbeta}

\newcommand{\var}{\mathrm{var}}
\newcommand{\cov}{\mathrm{cov}}

\newcommand{\Sig}{\mathbf{\Sigma}}

\newcommand{\tr}{\mathrm{tr}}
\newcommand{\diag}{\mathrm{diag}}

\newcommand{\supp}{\mathrm{supp}}

\newcommand{\argmax}[1]{\underset{#1}{\operatorname{arg}\,\operatorname{max}}\;}
\newcommand{\argmin}[1]{\underset{#1}{\operatorname{arg}\,\operatorname{min}}\;}

\newcommand{\lammin}{\lambda_{\min}}

\newcommand{\aic}{\mathrm{AIC}}

\newcommand{\gbicp}{$\mbox{GBIC}_p$ }
\newcommand{\gbicpl}{$\mbox{GBIC}_p$-L }
\newcommand{\bbetazero}{\mbox{\boldmath $\beta$}_{n,0}}

\def\boxit#1{\vbox{\hrule\hbox{\vrule\kern6pt\vbox{\kern6pt#1\kern6pt}\kern6pt\vrule}\hrule}}

\def\t{^T}
\def\toD{\overset{\mathscr{D}}{\longrightarrow}}

\begin{document}

\title{Model Selection in High-Dimensional Misspecified Models
\thanks{This work was partially supported by NSF CAREER Award DMS-0955316 and Grants DMS-0806030 and DMS-1308566.}
}
\author[1]{Pallavi Basu}
\author[2]{Yang Feng}
\author[1]{Jinchi Lv}
\affil[1]{
University of Southern California}
\affil[2]{
Columbia University}

\date{}
\maketitle
\begin{abstract}
Model selection is indispensable to high-dimensional sparse modeling in selecting the best set of covariates among a sequence of candidate models. Most existing work assumes implicitly that the model is correctly specified or of fixed dimensions. Yet model misspecification and high dimensionality are common in real applications. In this paper, we investigate two classical Kullback-Leibler divergence and Bayesian principles of model selection in the setting of high-dimensional misspecified models. Asymptotic expansions of these principles reveal that the effect of model misspecification is crucial and should be taken into account, leading to the generalized AIC and generalized BIC in high dimensions. With a natural choice of prior probabilities, we suggest the generalized BIC with prior probability which involves a logarithmic factor of the dimensionality in penalizing model complexity. We further establish the consistency of the covariance contrast matrix estimator in a general setting. Our results and new method are supported by numerical studies.


\end{abstract}

\section{Introduction}\label{sec::intro}
With rapid advances of modern technology, high-throughput data sets of
unprecedented size, such as genetic and proteomic data, fMRI and functional data, and panel data in economics and finance, are frequently encountered in many contemporary applications. In these applications, the dimensionality $p$ can be comparable to or even much larger than the sample size $n$. A key assumption that often makes large-scale inference feasible is the sparsity of signals, meaning that only a small fraction of covariates contribute to the response when $p$ is large compared to $n$. High-dimensional modeling with dimensionality reduction and feature selection plays an important role in these problems. A sparse modeling procedure typically produces a sequence of candidate models, each involving a possibly different subset of covariates. An important question is how to compare different models in high dimensions when models are possibly misspecified.

The problem of model selection has a long history with numerous contributions by many researchers. Among others, well-known model selection criteria are the AIC (Akaike, 1973 and 1974) and BIC (Schwarz, 1978), where the former is based on the Kullback-Leibler (KL) divergence principle of model selection and the latter is originated from the Bayesian principle. A great deal of work has been devoted to understanding and extending these methods. See, for example, Bozdogan (1987), Foster and George (1994), Konishi and Kitagawa (1996), Ing (2007), Chen and Chen (2008), Chen and Chan (2011), Ing and Lai (2011), Liu and Yang (2011), and Chang et al. (2014) in different model settings. The connections between the AIC and cross-validation have been investigated in Stone (1977), Hall (1990), and Peng et al. (2013) in various contexts. Model selection criteria such as AIC and BIC are frequently used for tuning parameter selection in regularization methods. For instance, mode selection in the context of penalized likelihood methods has been studied in Fan and Li (2001), Wang et al. (2007), Wang et al. (2009), Zhang et al. (2010), and Fan and Tang (2013). In particular, Fan and Tang (2013) showed that classical information criteria such as AIC and BIC can be inconsistent for model selection when the dimensionality $p$ grows very fast relative to sample size $n$.

Most existing work on model selection usually makes an implicit assumption that the model under study is correctly specified or of fixed dimensions. For example, White (1982) laid out a general theory of maximum likelihood estimation in misspecified models for the case of fixed dimensionality and independent and identically distributed (i.i.d.) observations. Recently, Lv and Liu (2014) investigated the problem of model selection with model misspecification and derived asymptotic expansions of both KL divergence and Bayesian principles in misspecified generalized linear models, leading to the generalized AIC and generalized BIC, for the case of fixed dimensionality. A specific form of prior probabilities motivated by the KL divergence principle leads to the generalized BIC with prior probability ($\mbox{GBIC}_p$-L\footnote{Here we use this notation to emphasize that the criterion is for the low-dimensional case, while reserving the original notation $\mbox{GBIC}_p$ in Lv and Liu (2014) for the high-dimensional counterpart.}).
Yet model misspecification and high dimensionality are both common in real applications. Thus a natural and important question is how to characterize the impact of model misspecification on model selection in high dimensions. We intend to provide some answer to this question in this paper. Our analysis enables us to suggest the generalized BIC with prior probability ($\mbox{GBIC}_p$) that involves a logarithmic factor of the dimensionality in penalizing model complexity.

To gain some insights into the challenges of the aforementioned problem, let us consider a motivating example. Assume that the response $Y$ depends on the covariate vector $(X_1, \cdots, X_p)\t$ through the functional form
\begin{align}\label{eq::multiple-index}
Y = f(X_1)+f(X_2-X_3)+f(X_4-X_5)+ \varepsilon,
\end{align}
where $f(x)=x^3/(x^2+1)$ and the remaining setting is as specified in Section \ref{sec::simu-multiple-index}. Consider sample size $n=100$ and vary  dimensionality $p$ from 200 to 3200. Without prior knowledge about the true model structure, we take the linear regression model
\begin{equation} \label{eq::linear-model}
\by = \bX \bbeta + \bveps
\end{equation}
as the working model, with the same notation therein, and apply some information criteria to hopefully recover the oracle working model consisting of the first five covariates. When $p=200$, the traditional AIC and BIC, which ignore model misspecification, tend to select a model with size larger than five. As expected, $\mbox{GBIC}_p$-L works reasonably well by selecting the oracle working model half of the time. However, when $p$ is increased to $3200$, these methods fail to select such a model with significant probability and the prediction performance of the selected models deteriorates. This motivates us to study the problem of model selection in high-dimensional misspecified models. In contrast, our newly suggested $\mbox{GBIC}_p$ can recover the oracle working model with significant probability in this challenging scenario.

The main contributions of our paper are threefold. First, we establish a systematic theory of model selection with model misspecification in high dimensions. The asymptotic expansions for different model selection principles involve delicate and challenging technical analysis. Second, our work provides rigorous theoretical justification of the covariance contrast matrix estimator that incorporates the effect of model misspecification and is crucial for practical implementation. Such an estimator is shown to be consistent in the general setting of high-dimensional misspecified models. Third, we suggest the use of a new prior in the expansion for $\mbox{GBIC}_p$ involving the $\log p$ term. This criterion has connections to the model selection criteria in Chen and Chen (2008) and Fan and Tang (2013) with the $\log p$ factor for the case of correctly specified models.

The rest of the paper is organized as follows. Section \ref{sec::definition} introduces the setup for model misspecification. We present some key asymptotic properties of the quasi-maximum likelihood estimator and provide asymptotic expansions of KL divergence and Bayesian model selection principles in high dimensions in Section \ref{sec::results}. Section \ref{sec::numerical} demonstrates the performance of different model selection criteria in high-dimensional misspecified models through several simulation  and real data examples. We provide some discussions of our results and possible extensions in Section \ref{sec::discussions}. The proofs of some main results are relegated to the Appendix. Additional technical proofs and numerical results are provided in the Supplementary Material.

\section{Model misspecification} \label{sec::definition}

Assume that conditional on the covariates, the $n$-dimensional random response vector $\bY = (Y_1, \cdots, Y_n)\t$ has a true unknown distribution $G_n$ with density function
\begin{equation} \label{001}
g_n(\by) = \prod_{i = 1}^n g_{n, i}(y_i),
\end{equation}
where $\by = (y_1, \cdots, y_n)\t$. Model (\ref{001}) entails that all components of $\bY$ are independent but not necessarily identically distributed. Consider a set of $d$ covariates out of all $p$ available covariates, where $p$ can be much larger than $n$. Denote by $\bX$ the corresponding $n \times d$ deterministic design matrix. To simplify the technical presentation, we focus on the case of deterministic design. In practice, one chooses a family of working models to fit the data. Model misspecification generally occurs when the family of distributions is misspecified or some true covariates are missed.

Since the true model $G_n$ is unknown, we choose a family of generalized linear models (GLMs) $F_n(\cdot, \bbeta) = F_n(\bz; \bX, \bbeta)$ with a canonical link as
our working models, each of which has density function
\begin{equation} \label{003}
f_n(\bz, \bbeta) d\mu_0(\bz)= \prod_{i = 1}^n f_0(z_i, \theta_i)
d\mu_0(z_i) \equiv \prod_{i = 1}^n \exp \left[z_i \theta_i  -
  b(\theta_i)\right] d\mu(z_i),
\end{equation}
where $\bz = (z_1, \cdots, z_n)\t$, $ \btheta = (\theta_1, \cdots,
\theta_n)\t = \bX \bbeta $
with $\bbeta \in \mR^d$, $b(\theta)$ is a smooth convex function, $\mu_0$ is
the Lebesgue measure, and $\mu$ is some fixed measure on
$\mR $. Assume that $b''(\theta)$ is continuous and bounded away from $0$, $\bX$ is of full column rank $d$, and $EY_i^2$ are bounded. Clearly $\{f_0(z, \theta): \theta \in \mR \}$ is a family of distributions in the regular exponential family and may
not contain $g_{n, i}$'s.


To ease the presentation, define two vector-valued functions
$\bb(\btheta) = (b(\theta_1), \cdots, b(\theta_n))\t$ and
$\bmu(\btheta) = (b'(\theta_1), \cdots,
b'(\theta_n))\t$, and a  matrix-valued function $\Sig(\btheta) =
\diag\{b''(\theta_1), \cdots,
b''(\theta_n)\}$. For any
$n$-dimensional random vector $\bZ$ with distribution $F_n(\cdot, \bbeta)$
given by (\ref{003}), it holds that
$E \bZ = \bmu(\bX \bbeta)$ and $\cov(\bZ) =
\Sig(\bX \bbeta)$.
The density function (\ref{003}) can be rewritten as
\[
f_n(\bz, \bbeta) = \exp \left[\bz\t \bX \bbeta - \bone\t \bb(\bX
  \bbeta)\right] \prod_{i = 1}^n \frac{d\mu}{d\mu_0}(z_i),
\]
where $\frac{d\mu}{d\mu_0}$ denotes the Radon-Nikodym
derivative. Given the observations $\by$ and $\bX$, this gives the
quasi-log-likelihood function
\begin{equation} \label{002}
\ell_n(\by, \bbeta) = \log f_n(\by, \bbeta) = \by\t \bX
  \bbeta - \bone\t \bb(\bX \bbeta) + \sum_{i = 1}^n \log
\frac{d\mu}{d\mu_0}(y_i).
\end{equation}
The quasi-maximum likelihood estimator (QMLE) of the
$d$-dimensional parameter vector $\bbeta$ is
defined as
\begin{equation} \label{005}
\hbbeta_n = \mbox{arg}\max_{\bbeta \in \mR^d} \ell_n(\by, \bbeta),
\end{equation}
which is the solution to the
score equation
$\bPsi_n(\bbeta) = \partial \ell_n(\by, \bbeta)/\partial \bbeta = \bX\t [\by - \bmu(\bX \bbeta)] = \bzero$.
This equation becomes the normal equation $\bX\t \by = \bX\t \bX \bbeta$ in
the linear regression model.

The KL divergence (Kullback and Leibler, 1951) of the model $F_n(\cdot,\bbeta)$ from the true model $G_n$ can be written as $I(g_n;f_n(\cdot, \bbeta)) = E \log g_n(\bY) -E\ell_n(\bY,\bbeta)$. The best working model that is closest to the true model under the KL divergence has parameter vector $\bbetazero=\arg\min_{\bbeta\in\mR^d}I(g_n;f_n(\cdot,\bbeta))$, which solves the equation
\begin{align}\label{eq:normal:equation}
\bX\t \left[E\bY -
\bmu(\bX \bbeta)\right] = \bzero.
\end{align}
We introduce two matrices that play a key role in model selection with model misspecification. Define
\begin{equation} \label{089}
\cov\left[\bPsi_n(\bbeta_{n, 0})\right] = \cov\left(\bX\t
\bY\right) = \bX\t \cov(\bY) \bX = \bB_n
\end{equation}
with $\cov(\bY)  = \diag\{\var(Y_1), \cdots, \var(Y_n)\}$ by the independence assumption,
\begin{equation} \label{090}
\frac{\partial^2 I(g_n; f_n(\cdot, \bbeta))}{\partial \bbeta^2}
= -\frac{\partial^2 \ell_n(\by, \bbeta)}{\partial \bbeta^2}
= \bX\t \Sig(\bX \bbeta) \bX = \bA_n(\bbeta),
\end{equation}
and $\bA_n = \bA_n(\bbeta_{n, 0})$. Observe that $\bA_n$ and $\bB_n$ are the covariance matrices of $\bX\t \bY$ under the best misspecified GLM $F_n(\cdot, \bbetazero)$ and the true model $G_n$, respectively.

\section{High-dimensional model selection in misspecified models} \label{sec::results}

We now present the asymptotic expansions of both KL divergence and Bayesian model selection principles in high-dimensional misspecified GLMs.

\subsection{Technical conditions and asymptotic properties of QMLE in high dimensions}



We list a few technical conditions required to prove the asymptotic properties of QMLE with diverging dimensionality. Denote by $\|\cdot\|_2$ the vector $L_2$-norm and the matrix operator norm.

\begin{assumption}\label{cond1}
There exists some constant $H > 0$ such that for each $1 \leq i \leq n$, $P(|q_i|> t)\leq H \exp(-t^2/H)$ for any $t\geq0$, where $(q_1,\cdots,q_n)^T=\cov(\bY)^{-1/2}(\bY-E\bY)$.
\end{assumption}

\begin{assumption}\label{condition:min_ev_B_n and min_ev_V}
There exist positive constants $c_1$, $c_0 > 8c_1^{-2}H$, and $r < 1/4$ such that for sufficiently large $n$, $\min_{\bbeta \in N_n(\delta_n)} \lambda_{\min} \left\{\bV_n(\bbeta)\right\} > c_1n^{-r}$ and $\lammin(\bB_n) \gg d\delta_n^2$, where $\delta_n = n^{r} (c_0 \log n)^{1/2}$, $N_n(\delta_n) = \{\bbeta \in \mR^d: \|(n^{-1} \bB_n)^{1/2} (\bbeta - \bbeta_{n, 0})\|_2 \leq (n/d)^{-1/2} \delta_n\}$, and $\bV_n(\bbeta) = \bB_n^{-1/2} \bA_n(\bbeta) \bB_n^{-1/2}$. Moreover, $d = o\{n^{(1-4r)/3} (\log n)^{-2/3}\}$.
\end{assumption}

\begin{assumption}\label{condition:Yi_deviation and sub-gaussian}
Assume $\sum_{i = 1}^n (\bx_i\t \bB_n^{-1}\bx_i)^{3/2} = o(1)$ and $\max_{1\leq i\leq n} E |Y_i - E Y_i|^3 = O(1)$.
\end{assumption}

\begin{assumption}\label{condition:widetilde_V_n}
Assume $$\max_{\bbeta_1, \cdots, \bbeta_d \in N_n(\delta_n)} \|\widetilde{\bV}_n(\bbeta_1, \cdots, \bbeta_d) - \bV_n\|_2 =  O(dn^{-1/2}\delta_n),$$ where $\bV_n = \bV_n(\bbeta_{n, 0}) = \bB_n^{-1/2} \bA_n \bB_n^{-1/2}$ and $\widetilde{\bV}_n(\bbeta_1, \cdots, \bbeta_d) = \bB_n^{-1/2} \widetilde{\bA}_n(\bbeta_1,$ $\cdots, \bbeta_d) \bB_n^{-1/2}$ with $\widetilde{\bA}_n(\bbeta_1, \cdots,\bbeta_d)$ a $d \times d$ matrix with $j$th row the corresponding row of $\bA_n(\bbeta_j)$ for each $1 \leq j \leq d$. Moreover, $\lambda_{\max}(\bV_n)$ is a polynomial order of $n$.
\end{assumption}

Conditions \ref{cond1} and \ref{condition:min_ev_B_n and min_ev_V} are some basic assumptions for establishing the consistency of the QMLE $\widehat\bbeta_n$ in Theorem \ref{Thm5}. In particular, Condition \ref{cond1} assumes that the standardized response has sub-Gaussian distribution which facilitates the derivation of the deviation probability bound. Conditions \ref{condition:min_ev_B_n and min_ev_V}--\ref{condition:widetilde_V_n} are similar to those in Lv and Liu (2014), except for some major differences due to the  high-dimensional setting. In particular, Condition \ref{condition:min_ev_B_n and min_ev_V} allows the minimum eigenvalue of $\bV_n(\bbeta)$ to converge to zero at a certain rate as $n$ increases in a  neighborhood $N_n(\delta_n)$ of $\bbeta_{n,0}$. Such a neighborhood is wider compared to that for the case of fixed dimensionality. The dimensionality  $d$ of the QMLE is allowed to diverge with $n$. Conditions \ref{condition:Yi_deviation and sub-gaussian} and \ref{condition:widetilde_V_n} are imposed to establish the asymptotic normality of $\widehat\bbeta_n$.



\begin{theorem} \emph{(Consistency of QMLE)}. \label{Thm5}
Under Conditions \ref{cond1}--\ref{condition:min_ev_B_n and min_ev_V}, the QMLE $\hbbeta_n$ satisfies $\hbbeta_n-\bbeta_{n,0}=o_P(1)$ and further $\hbbeta_n \in N_n(\delta_n)$ with probability $1-O(n^{-\alpha})$ for some large positive constant $\alpha$.
\end{theorem}

\begin{theorem} \emph{(Asymptotic normality)}. \label{Thm6}
Under Conditions \ref{cond1}--\ref{condition:widetilde_V_n}, the QMLE
$\hbbeta_n$ satisfies
\[ \bD_n \bC_n (\hbbeta_n - \bbeta_{n, 0}) \toD N(\bzero, I_m), \]
where $\bC_n = \bB_n^{-1/2} \bA_n$ and $\bD_n$ is any $m \times d$ matrix such that $\bD_n \bD_n^T = I_m$.
\end{theorem}

Theorems \ref{Thm5} and \ref{Thm6} establish the consistency and asymptotic normality of the QMLE in high-dimensional misspecified GLM. These results provide the theoretical foundation for the technical analyses in Sections \ref{Sec2.5}--\ref{Sec2.3}. The asymptotic theory of the QMLE reduces to that of the maximum likelihood estimator (MLE) when the model is correctly specified. Our results extend those in Lv and Liu (2014) for the case of fixed dimensionality. We next introduce a few additional conditions for deriving the asymptotic expansions of the two model selection principles.

\begin{assumption}\label{condition:region}
There exists some constant $\alpha_1$ with $0 < \alpha_1 < \alpha/2 - 1$ such that $b''(\cdot) = O(n^{\alpha_1})$ and for sufficiently large $n$, $N_n(\delta_n) \subset M_n (\alpha_1) = \{\bbeta \in \mR^d: \|\bX \bbeta \|_{\infty} \leq \alpha_1 \log n\}$, where constant $\alpha$ is given in Theorem \ref{Thm5}.
\end{assumption}

\begin{assumption}\label{condition:prior density and rho_n}
Assume that 
$\pi(h(\bbeta)) = \frac{d\mu_{\mathfrak{M}}}{d\mu_0} (h(\bbeta))$ satisfies
\begin{equation}\label{condonpi}
\inf_{\bbeta \in N_n(2\delta_n)} \pi (h(\bbeta)) \geq c_2 \text{  and  } \sup_{\bbeta \in \mathbb{R}^d} \pi (h(\bbeta)) \leq c_3
\end{equation}
with $c_2, c_3 > 0$ some constants, and
$\rho_n(\delta_n) = \max_{\bbeta \in N_n(2\delta_n)} \max \{ | \lambda_{\min}(\bV_n(\bbeta) \\- \bV_n) |, | \lambda_{\max} (\bV_n(\bbeta) - \bV_n) | \} = o\{n^{-(1-r)/3}\}$.
\end{assumption}

\begin{assumption}\label{condition:lipschitz}
Assume that $n^{-1} \bA_n(\bbeta)$, $n^{-1}\bX^T\diag\{|\bmu(\bX\bbeta) - \bmu(\bX\bbeta_{n,0})|\}\\ \bX$, and $n^{-1}\bX^T\diag\{[\bmu(\bX\bbeta) - \bmu(\bX\bbeta_{n,0})] \circ [\bmu(\bX\bbeta) - \bmu(\bX\bbeta_{n,0})] \}\bX$ are Lipschitz (in operator norm) with constant $L > 0$ in $N_n(\delta_n)$, and $\| \bX \|_{\infty} = O(n^{\alpha_2})$ with constant $0 \leq \alpha_2 < r$, where $\circ$ represents the Hadamard (componentwise) product and $\|\cdot\|_\infty$ denotes the entrywise matrix $L_{\infty}$-norm.
\end{assumption}

\begin{assumption}\label{condition:bias_variance}
Assume $\sum_{i=1}^n \{[EY_i - (\bmu(\bX\bbeta_{n,0}))_i]^2/\var(Y_i)\}^2 = O(n^{\alpha_3})$ with some constant $0 \leq \alpha_3 \leq 4(r - \alpha_2)$.
\end{assumption}


The first part of Condition \ref{condition:region} holds naturally for linear and logistic regression models, and is introduced to accommodate the case of Poisson regression. The second part of Condition \ref{condition:region} is a mild assumption ensuring that the restricted QMLE coincides with its unrestricted version with significant probability, which is key to the asymptotic expansion of the KL divergence principle in high dimensions in Theorem \ref{Thm3}. It is worth mentioning that the set $M_n (\alpha_1)$ grows with $n$, while the neighborhood $N_n(\delta_n)$ is asymptotically shrinking.

Condition \ref{condition:prior density and rho_n} is similar to the one in Lv and Liu (2014), except that we need to specify the rate at which $\rho_n(\delta_n)$ converges to zero. Condition \ref{condition:lipschitz} requires the Lipschitz property for those matrix-valued functions. The bound on the entry-wise matrix $L_\infty$-norm of the design matrix is mild. Condition \ref{condition:bias_variance} is a sensible assumption bounding the effect of the model bias. In particular, Conditions \ref{condition:lipschitz} and \ref{condition:bias_variance} are introduced only for proving the consistency of the covariance contrast matrix in the general setting in Theorem \ref{Thm4}.

\subsection{Generalized AIC in misspecified models} \label{Sec2.5}
Given a sequence of subsets $\left\{\mathfrak{M}_m: m
= 1, \cdots, M\right\}$ of the full model $\{1, \cdots, p\}$, we can
construct a sequence of QMLE's $\{\hbbeta_{n, m}: m = 1, \cdots, M\}$
by fitting the GLM (\ref{003}). A natural
question is how to compare those fitted models. The QMLEs $\{\hbbeta_{n, m}: m = 1, \cdots, M\}$ become the MLEs when the model is correctly specified.

Akaike's principle of model
selection is choosing
the model $\mathfrak{M}_{m_0}$ that minimizes the KL
divergence $I(g_n; f_n(\cdot, \hbbeta_{n, m}))$ of the fitted
model  $F_n(\cdot, \hbbeta_{n, m})$ from the true model $G_n$, that is,
\begin{equation} \label{060}
m_0 = \argmin{m \in \{1, \cdots, M\}} I(g_n; f_n(\cdot, \hbbeta_{n, m})),
\end{equation}
where
\begin{equation} \label{012}
I(g_n; f_n(\cdot, \hbbeta_{n, m})) = E \log g_n(\widetilde{\bY}) - \eta_n(\hbbeta_{n, m})
\end{equation}
with $\eta_n(\bbeta) = E \ell_n(\widetilde{\bY}, \bbeta)$ and
$\widetilde{\bY}$ an independent
copy of $\bY$. Thus
\[
m_0 = \argmax{m \in \{1, \cdots, M\}} \eta_n(\hbbeta_{n, m}) = \argmax{m \in \{1, \cdots, M\}} E_{\widetilde{\bY}}
\ell_n(\widetilde{\bY}, \hbbeta_{n, m}),
\]
which shows that Akaike's principle of model selection is equivalent to choosing the
model  $\mathfrak{M}_{m_0}$ that maximizes the expected
log-likelihood with the expectation taken with respect to an independent copy of $\bY$. Using the asymptotic theory of MLE, Akaike (1973)
showed that for the case of i.i.d. observations,
$\eta_n(\hbbeta_n)$ can be asymptotically
expanded as $\ell_n(\by, \hbbeta_n) -
|\mathfrak{M}|$, which leads to the seminal AIC for comparing
competing models:
\begin{equation} \label{062}
\aic(\by, \mathfrak{M}) = -2 \ell_n(\by, \hbbeta_n) + 2 |\mathfrak{M}|.
\end{equation}
For simplicity, we drop the last term in (\ref{002}) which
does not depend on $\bbeta$, and redefine the quasi-log-likelihood as
$\ell_n(\by, \bbeta) = \by\t \bX \bbeta - \bone\t \bb(\bX \bbeta)$ hereafter.

\begin{theorem} \label{Thm3}
Under Conditions \ref{cond1}--\ref{condition:region}, we have with probability tending to one,
\begin{align}
  E\eta_n(\widehat\bbeta_n)=E\ell_n(\by,\widehat\bbeta_n)-\tr(\bH_n)+o\{\tr(\bH_n)\},
\end{align}
where $\bH_n=\bA_n^{-1} \bB_n$.
\end{theorem}

Theorem \ref{Thm3} generalizes the corresponding result in Lv and Liu (2014) to high dimensions. However, we would like to point out that our new technical analysis differs substantially from theirs due to the challenges of diverging dimensionality. The asymptotic expansion in Theorem \ref{Thm3}  enables us to introduce the generalized AIC (GAIC) as follows.

\begin{definition}\label{def::GAIC}
We define $GAIC$ of model $\mathfrak{M}$ as
\begin{align}\label{def::GAIC_p}
GAIC(\by,\mathfrak{M};F_n)= -2 \ell_n(\by, \hbbeta_n) + 2\tr(\widehat\bH_n),
\end{align}
where $\widehat\bH_n$ is a consistent estimator of $\bH_n$ specified in Section \ref{sec::est-H-n}.
\end{definition}
When the model is correctly specified, it holds that  $\tr(\widehat\bH_n)\approx \tr(\bI_d)=|\mathfrak{M}|$, under which GAIC reduces to AIC asymptotically. We demonstrate in the simulation studies that GAIC can improve over the original AIC substantially in the presence of model misspecification.

%
%

\subsection{Estimation of covariance contrast matrix} \label{sec::est-H-n}

From the asymptotic expansions for the GAIC, GBIC, and $\mbox{GBIC}_p$ (the latter two to be introduced in Section \ref{Sec2.3}), a common term is the covariance contrast matrix $\bH_n$, which characterizes the impact of model misspecification. Therefore, providing an accurate estimator for such a matrix $\bH_n$ is of vital importance in the application of these information criteria.

Consider the plug-in estimator $\widehat\bH_n=\widehat\bA_n^{-1}\widehat\bB_n$ with $\widehat\bA_n$ and $\widehat\bB_n$ defined as follows. Since the QMLE $\hbbeta_n$ provides a consistent estimator of $\bbeta_{n, 0}$ in the best misspecified GLM $F_n(\cdot, \bbeta_{n, 0})$, a natural estimate of matrix $\bA_n$ is given by
\begin{equation} \label{eq::A_n_hat}
\widehat{\bA}_n = \bA_n(\hbbeta_n) = \bX\t \Sig(\bX \hbbeta_n) \bX.
\end{equation}
When the model is correctly specified, the following simple estimator
\begin{equation} \label{eq::B_n_hat}
\widehat{\bB}_n = \bX\t \diag\left\{\left[\by - \bmu(\bX
  \hbbeta_n)\right] \circ \left[\by - \bmu(\bX
  \hbbeta_n)\right]\right\} \bX
\end{equation}
gives an asymptotically unbiased estimator of $\bB_n$.

\begin{theorem} (Consistency of estimator) \label{Thm4}
Assume that Conditions \ref{cond1}--\ref{condition:Yi_deviation and sub-gaussian} and \ref{condition:lipschitz}--\ref{condition:bias_variance} hold, the eigenvalues of $n^{-1}\bA_n$ and $n^{-1}\bB_n$ are bounded away from 0 and $\infty$, and $d = o\{n^{(1-4r)/4}\}$. Then the plug-in estimator $\widehat\bH_n$ satisfies $\tr(\hH_n) = \tr(\bH_n) + o_P(1)$ and $\log |\hH_n| = \log |\bH_n| + o_P(1)$.
\end{theorem}

Theorem \ref{Thm4} improves the result in Lv and Liu (2014) in two important aspects. First, the consistency of the covariance contrast matrix estimator was previously justified in Lv and Liu (2014) for the case of correctly
specified model. Our new result shows that the simple plug-in estimator $\widehat\bH_n$ still enjoys consistency in the general setting of model misspecification. Second, the result in Theorem \ref{Thm4} holds for the case of diverging dimensionality. These theoretical guarantees are crucial to the practical implementation of those information criteria. Our numerical studies reveal that such an estimate works well in a variety of model misspecification settings.

\subsection{Generalized BIC in misspecified models} \label{Sec2.3}
Given a set of competing models $\left\{\mathfrak{M}_m: m = 1,
\cdots, M\right\}$, a popular Bayesian model selection procedure is to first put nonzero
prior probability $\alpha_{\mathfrak{M}_m}$ on each  model
$\mathfrak{M}_m$, and then choose a prior distribution
$\mu_{\mathfrak{M}_m}$ for the parameter vector in the corresponding model. Assume
that the density function of $\mu_{\mathfrak{M}_m}$  is bounded in
$\mR^{\mathfrak{M}_m}
 = \mR^{d_m}$ with $d_m=|\mathfrak{M}_m|$  and locally bounded away from
 zero throughout the domain. The Bayesian principle of model selection  is
to choose the most probable model {\it a
  posteriori}, that is, choose model  $\mathfrak{M}_{m_0}$ such that
\begin{equation} \label{021}
m_0 = \arg \max_{m \in \{1, \cdots, M\}} S(\by, \mathfrak{M}_m; F_n),
\end{equation}
where the log-marginal-likelihood is
\begin{equation}\label{083}
 S(\by, \mathfrak{M}_m; F_n) = \log \int \alpha_{\mathfrak{M}_m}
\exp\left[\ell_n(\by, \bbeta)\right] d\mu_{\mathfrak{M}_m}(\bbeta)
\end{equation}
with the log-likelihood $\ell_n(\by, \bbeta)$ as in
(\ref{002}) and the integral over $\mR^{d_m}$.

To ease the presentation, for any $\bbeta\in \mR^d$ we define a quantity
\begin{align}\label{eq:ellnstar}
\ell_n^*(\by,\bbeta)=\ell_n(\by,\bbeta)-\ell_n(\by,\widehat\bbeta_n),
\end{align}
which is the deviation of the quasi-log-likelihood from its maximum. Then from (\ref{083}) and (\ref{eq:ellnstar}), we have
\begin{equation} \label{184}
S(\by, \mathfrak{M}_m; F_n) = \ell_n(\by,\widehat\bbeta_n) + \log E_{\mu_{\mathfrak{M}_m}} [ U_n(\bbeta)^n ] + \log \alpha_{\mathfrak{M}_m},
\end{equation}
where $U_n(\bbeta) = \exp[n^{-1} \ell^*_n(\by, \bbeta)]$.

\begin{theorem} \label{Thm1}
Under Conditions \ref{cond1}--\ref{condition:Yi_deviation and sub-gaussian} and \ref{condition:prior density and rho_n}, we have with probability tending to one,
\begin{align}
S(\by, \mathfrak{M}; F_n) & = \ell_n(\by, \hbbeta_n) -\frac{\log n}{2}
|\mathfrak{M}| + \frac{1}{2} \log |\bH_n| + \log
\alpha_\mathfrak{M} \\
& \quad + \frac{\log (2\pi)}{2}|\mathfrak{M}|+\log c_n+o(1), \nonumber
\end{align}
where $\bH_n=\bA_n^{-1} \bB_n$ and $c_n \in [c_2,c_3]$.
\end{theorem}

The  asymptotic expansion of the Bayes factor in Theorem \ref{Thm1} leads us to introduce the generalized BIC (GBIC) as follows.

\begin{definition}
We define GBIC of model $\mathfrak{M}$ as
\begin{align}\label{def::GBIC}
GBIC(\by,\mathfrak{M};F_n)= -2 \ell_n(\by, \hbbeta_n) +(\log n) |\mathfrak{M}|-\log|\widehat\bH_n|,
 \end{align}
where $\widehat\bH_n$ is a consistent estimator of $\bH_n$.
\end{definition}

It is clear from \eqref{def::GBIC} that GBIC contains an extra term compared to BIC that replaces the factor $2$ with $\log n$ in penalizing model complexity in (\ref{062}). This additional term reflects the effect of model misspecification. When the model is correctly specified, GBIC reduces to BIC asymptotically.

The choice of the prior probabilities $\alpha_{\mathfrak{M}_{m}}$ is important in high dimensions. Lv and Liu (2014) suggested prior probability $\alpha_{\mathfrak{M}_{m}} \propto e^{-D_{m}}$ for each candidate model $\mathfrak{M}_m $, where the quantity $D_m$ is defined as
\begin{equation}
D_m = E \left[I(g_n; f_n(\cdot,\hbbeta_{n,m})) - I(g_n; f_n(\cdot,\bbeta_{n,m,0}))  \right]
\end{equation}
and the subscript $m$ indicates a particular candidate model. The motivation is that the further the QMLE $\hbbeta_{n,m}$ is away from the best misspecified GLM $F_n(\cdot, \bbeta_{n,m,0})$, the lower prior we assign to that model. In the high-dimensional setting when $p$ can be much larger than $n$, it is sensible to take into account the complexity of the space of all possible sparse models with the same size as $\mathfrak{M}_m$. This observation motivates us to consider a new prior of the form
\begin{equation} \label{prior}
\alpha_{\mathfrak{M}_{m}} \propto {{p\choose d}}^{-1} e^{-D_m}
\end{equation}
with $d=|\mathfrak{M}_{m}|$. Such a complexity factor has been exploited in the extended BIC (EBIC) in Chen and Chen (2008), who showed that using the term ${{p\choose d}}^{-\gamma}$ with some constant $0 < \gamma \leq 1$, the EBIC  can be model selection consistent for $p=O(n^{\kappa})$ with some positive constant $\kappa$ satisfying $1-(2\kappa)^{-1} < \gamma$.

Under the assumption of $d = o(p)$, an application of Stirling's formula shows that up to an additive constant, it holds that $ \log \alpha_{\mathfrak{M}_{m}} \approx -D_m-d\log p - d + d \log d$.
Thus for the prior defined in \eqref{prior}, we have an additional term $-(\log p+1-\log d)|\mathfrak{M}|$ in the asymptotic expansion for GBIC. When $p$ is of order $n^{\kappa}$ with some constant $\kappa > 0$, this new term is of the same order as $-(\log n) |\mathfrak{M}|$. When $\log p$ is of order $n^{\kappa}$ with some constant $0<\kappa<1$, the $\log p$ term dominates that involving $\log n$. Fan and Tang (2013) proposed a similar term $\log(\log n) \log p$ term to ameliorate the BIC for the case of correctly specified models with non-polynomially growing dimensionality $p$. The following theorem provides the asymptotic expansion of the Bayes factor with the particular choice of prior in (\ref{prior}).

\begin{theorem}  \label{Thm2}
Assume that Conditions \ref{cond1}--\ref{condition:prior density and rho_n} hold, $\alpha_{\mathfrak{M}_m}=C {{p\choose d}}^{-1}e^{-D_m}$ with $C > 0$ some normalization constant, and $d = o(p)$. Then we have with probability tending to one,
\begin{align}
S(\by, \mathfrak{M}; F_n) & = \ell_n(\by, \hbbeta_n) -(\log p^*)
|\mathfrak{M}| - \frac{1}{2}\tr(\bH_n) + \frac{1}{2} \log |\bH_n| \\
& \quad + \log (Cc_n)+o(1), \nonumber
\end{align}
where $\bH_n=\bA_n^{-1} \bB_n$, $p^* = \max\{n, p\}$, and $c_n \in [c_2,c_3]$.
\end{theorem}

Similarly to the GBIC, we now define a new information criterion, the  generalized BIC with prior probability ($\mbox{GBIC}_p$), based on Theorem \ref{Thm2}.

\begin{definition}
We define $GBIC_p$ of model $\mathfrak{M}$ as
\begin{align}\label{def::GBIC_p}
GBIC_p(\by,\mathfrak{M};F_n)= -2 \ell_n(\by, \hbbeta_n) +2 (\log p^*) |\mathfrak{M}|+\tr(\widehat\bH_n)-\log|\widehat\bH_n|,
 \end{align}
where $\widehat\bH_n$ is a consistent estimator of $\bH_n$.
\end{definition}

In correctly specified models, the term $\tr(\widehat\bH_n)-\log|\widehat\bH_n|$ is asymptotically close to $|\mathfrak{M}|$ when $\widehat\bH_n$ is a consistent estimator of $\bH_n=\bI_d$. Thus compared to BIC with factor $\log n$, the $\mbox{GBIC}_p$ contains a larger factor $\log p$ when $p$ grows non-polynomially with $n$. This leads to a heavier penalty on model complexity similarly as in Fan and Tang (2013). As pointed out in Lv and Liu (2014), the right hand side of (\ref{def::GBIC_p}) can be viewed as a sum of three terms: the goodness of fit, model complexity, and model misspecification. An important distinction with the low-dimensional counterpart of $\mbox{GBIC}_p$ is that our new criterion explicitly takes into account the dimensionality of the whole feature space.

\section{Numerical studies} \label{sec::numerical}
The asymptotic expansions of both KL divergence and Bayesian principles in Section \ref{sec::results} have enabled us to introduce the GAIC, GBIC, and $\mbox{GBIC}_p$ for model selection in high dimensions with model misspecification. We now investigate their performance in comparison to the information criteria AIC, BIC, and $\mbox{GBIC}_p$-L in high-dimensional misspecified models via simulation examples as well as two real data sets. For each simulation study, we set the number of repetitions to be 100 and examined the scenarios when the dimensionality grows ($p=200$, 400, 1600, and 3200).

\subsection{Simulation examples} \label{sec::simulation}

\subsubsection{Sparse linear regression with interaction and weak effects}\label{sec::simu-linear-interaction-weak}
The first model we consider is the following  high-dimensional linear regression model with interaction and weak effects
\begin{equation} \label{eq::simu-linear-interaction}
\by = \bX \bbeta + \bx_{p + 1} + \bveps,
\end{equation}
where $\bX = (\bx_1, \cdots, \bx_p)$ is an $n \times p$ design matrix, $\bx_{p+1} = \bx_1 \circ \bx_2$ is an interaction term which is the product of the first two covariates,  the rows of $\bX$ are sampled as i.i.d. copies from $N(\bzero, I_p)$,  and the error vector $\bveps\sim N(\bzero, \sigma^2 I_n)$. We set $\bbeta_0 = (1, -1.25, 0.75, -0.95, 1.5, 0.1,\\-0.1,0.1,-0.1,0.1, 0, \cdots, 0)\t$, $n = 100$, and $\sigma = 0.25$. Although the data was generated from model (\ref{eq::simu-linear-interaction}), we fit the linear regression model \eqref{eq::linear-model} without interaction, which is a typical example of model misspecfication. In view of (\ref{eq::simu-linear-interaction}), the true model involves only the first ten covariates in a nonlinear form. Since the other covariates are independent of those ten covariates, the oracle working model is $\supp(\bbeta_0) = \{1, \cdots, 10\}$ as argued in Lv and Liu (2014). Due to the high dimensionality,  it is computationally prohibitive to implement the best subset selection. Therefore, we first applied the regularization method SICA  (Lv and Fan, 2009) to build a sequence of sparse models and then selected the final model using a model selection criterion. In practice, one can apply any preferred variable selection procedure to obtain a sequence of candidate models.

In addition to comparing the models selected by different information criteria, we also considered the estimate based on the oracle working model $M_0=\{1, \cdots, 10\}$ as a benchmark and used both measures of prediction and variable selection. Denote by $\widehat M$ the selected model.
We split the oracle working model into the set of strong effects $M_{0,s}=\{1,\cdots,5\}$ and that of weak effects $M_{0,w}=\{6,\cdots,10\}$.  It is interesting to observe that all criteria tend to miss the entire set of weak effects $M_{0,w}$ due to their very low signal strength. Therefore, we focused on comparing the model selection performance in recovering the set of  strong effects $M_{0,s}$.

We report the strong effect consistent selection probability (the portion of simulations where $\widehat M=M_{0,s}$), the strong effect inclusion probability (the portion of simulations where $\widehat M\supset M_{0,s}$), and the prediction error $E (Y - \bx\t \hbbeta)^2$ with $\hbbeta$ an estimate and $(\bx\t, Y)$ an independent observation.
To evaluate the prediction performance of different criteria, we calculated the average prediction error on  an independent test sample of size 10,000. The results for prediction error and model selection performance are summarized in Table \ref{tb::simu-linear-inter-weak}. To save space, the number of false positives $|\widehat M \cap M_{0}^c|$ and the numbers of false negatives for strong effects $|\widehat M^c \cap M_{0,s}|$ and weak effects $|\widehat M^c \cap M_{0,w}|$, respectively, are reported in Table \ref{tb::simu-linear-inter-weak-FP-FN} in the Supplementary Material.

It is clear that as the dimensionality $p$ increases, the consistent selection probability tends to decrease and the prediction error tends to increase  for all information criteria. Generally speaking, GAIC improved over AIC, and GBIC, $\mbox{GBIC}_p$-L, and $\mbox{GBIC}_p$ performed better than BIC in terms of both prediction and variable selection. In particular, the model selected by our new information criterion $\mbox{GBIC}_p$ delivered the best performance with the smallest prediction error and highest strong effect consistent selection probability across all settings.

Meanwhile it is also interesting to see what results different model selection criteria lead to when the model is correctly specified. To this end, we regenerate the solution path based on the linear regression model with the interaction $\bx_{p+1}=\bx_1\circ \bx_2$ added. The same performance measures are calculated for this scenario with the results reported in Tables \ref{tb::simu-linear-inter-weak-correct} and \ref{tb::simu-linear-inter-weak-FP-FN-correct}, where the latter table is included in the Supplementary Material. A comparison of these results with those in Tables \ref{tb::simu-linear-inter-weak} and \ref{tb::simu-linear-inter-weak-FP-FN} gives several interesting observations. First, all model selection criteria have a better performance when the model is correctly specified in terms of both model selection and prediction. Second, it is worth noting that while all model selection criteria except AIC  work reasonably well for the correctly specified model, all but the newly proposed \gbicp have a very low consistent selection probability under both model misspecification and high dimensionality. Third, it is interesting to see that \gbicp outperforms the existing methods even under the correctly specified model in terms of consistent selection probability.

\begin{table}\caption{Simulation results for Example \ref{sec::simu-linear-interaction-weak} with all entries multiplied by 100 when the model is misspecified, with the oracle results based on both strong effects and weak effects.\label{tb::simu-linear-inter-weak}}
\begin{center}
{

\begin{tabular}{c|ccccccc}

\hline
\multicolumn{8}{c}{Strong effect consistent selection probability with inclusion probability }\\
\hline
 $p$ &  AIC & BIC & GAIC & GBIC & $\mbox{GBIC}_p$-L & $\mbox{GBIC}_p$&Oracle \\
\hline
200 &   0(99) &  29(99) &  21(99) &  32(99) &  67(98) &  73(98)&100(100) \\
400 &   0(100) &   9(100) &   8(100) &  19(100) &  54(100) &  76(100)&100(100) \\
1600 &   0(100) &   0(100) &   9(100) &   0(100) &  27(100) &  66(100)&100(100) \\
3200 &   0(100) &   0(100) &   4(100) &   0(100) &  16(100) &  64(100) &100(100)\\
\hline
\multicolumn{8}{c}{Median prediction error with robust standard deviation in parentheses 
}\\
\hline
200 & 164(35) & 130(13) & 130(10) & 128(12) & 125(8) & 125(8) & 121(7) \\
400 & 162(29) & 154(38) & 129(13) & 131(22) & 125(9) & 122(10) & 120(7) \\
1600 & 168(31) & 172(28) & 134(13) & 170(28) & 129(14) & 125(10) & 121(7) \\
3200 & 159(22) & 169(23) & 135(14) & 167(23) & 134(15) & 125(13) & 120(8) \\
\hline
\end{tabular}
}

\end{center}

\end{table}

\begin{table}\caption{Simulation results for Example \ref{sec::simu-linear-interaction-weak} with all entries multiplied by 100 when the model is correctly specified, with the oracle results based on both strong effects and weak effects.\label{tb::simu-linear-inter-weak-correct}}
\begin{center}
\begin{tabular}{c|ccccccc}

\hline
\multicolumn{8}{c}{Strong effect consistent selection probability  with inclusion probability }\\
\hline
 $p$ &  AIC & BIC & GAIC & GBIC & $\mbox{GBIC}_p$-L & $\mbox{GBIC}_p$&Oracle \\
\hline
200 &   2(100) &  82(100) &  81(99) &  82(100) &  87(100) &  91(100) &100(100)\\
400 &   8(100) &  76(100) &  76(100) &  84(100) &  90(100) &  94(100)&100(100) \\
1600 &  39(95) &  74(99) &  65(89) &  79(99) &  88(100) &  96(100)&100(100) \\
3200 &  64(94) &  84(98) &  72(88) &  84(98) &  94(100) &  95(100) &100(100)\\
\hline
\multicolumn{8}{c}{Median prediction error with robust standard deviation (RSD)  in parentheses 
}\\
\hline
200 & 13.6(1.9) & 11.2(1.0) & 11.2(1.0) & 11.2(1.0) & 11.6(1.3) & 11.7(1.2) & 7.0(0.4) \\
400 & 12.1(1.4) & 11.5(1.3) & 11.5(1.2) & 11.5(1.2) & 11.7(1.3) & 11.8(1.0) & 6.9(0.4) \\
1600 & 12.4(8.3) & 12.0(7.9) & 11.9(9.8) & 12.0(8.0) & 12.2(7.7) & 12.4(7.3) & 7.0(0.4) \\
3200 & 21.2(10.2) & 20.7(9.4) & 21.8(11.0) & 20.7(9.4) & 20.4(8.8) & 20.3(8.5) & 7.0(0.3) \\
\hline
\end{tabular}

\end{center}

\end{table}

\normalsize

\subsubsection{Multiple index model}\label{sec::simu-multiple-index}
We next consider another model misspecification setting that involves the multiple index model
\begin{align}\label{simu::multiple-index}
Y = f(\beta_1X_1)+f(\beta_2X_2+\beta_3X_3)+f(\beta_4X_4+\beta_5X_5)+\varepsilon,
\end{align}
where the response depends on the covariates only through the first five ones  but with nonlinear functions and $f(x)=x^3/(x^2+1)$. Here the design matrix $\bX = (\bx_1, \cdots, \bx_p)$ was generated as in Section \ref{sec::simu-linear-interaction-weak}. We set the true parameter vector $\bbeta_0 = (1, 1, -1, 1, -1, 0, \cdots, 0)\t$, $n = 100$, and $\sigma = 0.25$. Note that the oracle working model is $M_0=\supp(\bbeta_0) = \{1, \cdots, 5\}$ for this example. Although the data was generated from model (\ref{simu::multiple-index}), we fit the linear regression model \eqref{eq::linear-model}. The results are summarized in Tables \ref{tb::simu-multiple-index} and \ref{tb::simu-multiple-index-FP-FN} (the latter available in Supplementary Material). The consistent selection probability and inclusion probability are now calculated based on $M_0$.

In general, the conclusions are similar to those in Example \ref{sec::simu-linear-interaction-weak}. An interesting observation is the comparison between \gbicpl  and \gbicp in terms of model selection. While \gbicpl is comparable to \gbicp when the dimension is not large ($p=200$), the difference between these two methods increases as the dimensionality increases. In the case when $p=3200$, \gbicp has 77\% success probability of consistent selection, while all the other criteria have at most 5\% success probability.  This confirms the necessity of including the $\log p$ factor in the model selection criterion to take into account the high dimensionality, which is in line with the conclusion in Fan and Tang (2013) for the case of correctly specified models.

\begin{table}\caption{Simulation results for Example \ref{sec::simu-multiple-index} with all entries multiplied by 100.  \label{tb::simu-multiple-index}}
\begin{center}

\begin{tabular}{c|ccccccc}
\hline
\multicolumn{8}{c}{Consistent selection probability  with inclusion probability }\\
\hline
 $p$ &  AIC & BIC & GAIC & GBIC & $\mbox{GBIC}_p$-L & $\mbox{GBIC}_p$&Oracle \\
\hline
200 & 2(100) & 4(100) & 2(100) & 6(100) & 51(100) & 65(100) &100(100)\\
400 & 1(100) & 1(100) & 2(100) & 1(100) & 28(100) & 67(100) &100(100)\\
1600 & 0(100) & 0(100) & 3(100) & 0(100) & 5(100) & 63(100) &100(100)\\
3200 & 0(100) & 0(100) & 5(100) & 0(100) & 5(100) & 77(100) &100(100)\\
\hline
\multicolumn{8}{c}{Median prediction error with RSD in parentheses}\\
\hline
200 & 26(3) & 26(3) & 26(3) & 26(3) & 23(3) & 23(2) & 22(1) \\
400 & 28(3) & 28(3) & 27(3) & 28(3) & 25(4) & 23(2) & 22(1) \\
1600 & 31(3) & 31(3) & 30(4) & 31(3) & 30(4) & 23(4) & 22(1) \\
3200 & 31(4) & 31(4) & 30(3) & 31(4) & 30(3) & 23(2) & 22(1) \\
\hline
\end{tabular}

\end{center}

\end{table}

\subsubsection{Logistic regression with interaction}\label{sec::simu-logi-inter}
Our last simulation example is high-dimensional logistic regression with interaction. We simulated 100 data sets from the logistic regression model with interaction and an $n$-dimensional parameter vector
\begin{equation} \label{096}
\btheta = \bX \bbeta + 2 \bx_{p + 1} + 2\bx_{p+2},
\end{equation}
where $\bX = (\bx_1, \cdots, \bx_p)$ is an $n \times p$ design matrix, $\bx_{p+1} = \bx_1 \circ \bx_2$ and $\bx_{p+2} = \bx_3 \circ \bx_4$ are two interaction terms, and the rest is the same as in (\ref{eq::simu-linear-interaction}). For each data set, the $n$-dimensional response vector $\by$ was sampled from the Bernoulli distribution with success probability vector $[e^{\theta_1}/(1+e^{\theta_1}), \cdots, e^{\theta_n}/(1+e^{\theta_n})]\t$ with $\btheta = (\theta_1, \cdots, \theta_n)\t$ given in (\ref{096}). As in Section \ref{sec::simu-linear-interaction-weak}, we consider the case where all covariates are independent of each other. We chose $\bbeta_0 = (2.5, -1.9, 2.8, -2.2, 3, 0, \cdots, 0)\t$ and set sample size $n = 200$. Although the data was generated from the logistic regression model with parameter vector (\ref{096}), we fit the logistic regression model without the two interaction terms. This provides another example of misspecified models. As argued in Section \ref{sec::simu-linear-interaction-weak}, the oracle working model is $\supp(\bbeta_0) = \{1, \cdots, 5\}$ which corresponds to the logistic regression model  with the first five covariates.

Since the goal in logistic regression is usually classification, we replace the prediction error with the classification error rate.  Tables \ref{tb::simu-logi-inter} and \ref{tb::simu-logi-inter-FP-FN} (the latter available in Supplementary Material) show similar phenomenon as in Sections \ref{sec::simu-linear-interaction-weak} and \ref{sec::simu-multiple-index}. Again $\mbox{GBIC}_p$ outperformed all other model selection criteria with greater advantage for the high-dimensional case (e.g., $p=3200$).

\begin{table}\caption{Simulation results for Example \ref{sec::simu-logi-inter} with all entries multiplied by 100.  }
\label{tb::simu-logi-inter}

\begin{center}

\begin{tabular}{c|ccccccc}
\hline
\multicolumn{8}{c}{Consistent selection probability  with inclusion probability }\\
\hline
  $p$&  AIC & BIC & GAIC & GBIC & $\mbox{GBIC}_p$-L & $\mbox{GBIC}_p$&Oracle \\
\hline
200 & 0(99) & 32(94) & 1(99) & 39(94) & 49(91) & 49(91)&100(100) \\
400& 0(99) & 19(97) & 0(99) & 36(93) & 50(92) & 55(92) &100(100)\\
1600& 0(96) & 0(96) & 0(94) & 21(90) & 35(88) & 47(81)&100(100) \\
3200 & 0(95) & 0(95) & 0(96) & 10(90) & 21(86) & 41(72)&100(100) \\
\hline
\multicolumn{8}{c}{Median classification error rate with RSD in parentheses}\\
\hline
200 & 22(3) & 15(2) & 16(2) & 15(1) & 14(1) & 14(1) & 14(1) \\
400 & 21(3) & 16(5) & 17(2) & 15(1) & 15(1) & 15(1) & 13(1) \\
1600 & 21(2) & 21(2) & 18(1) & 16(3) & 15(1) & 16(2) & 14(1) \\
3200 & 22(2) & 21(2) & 19(2) & 18(3) & 15(2) & 15(2) & 13(1) \\
\hline
\end{tabular}

\end{center}

\end{table}

\subsection{Real data examples} \label{Sec4.2}
We finally consider two gene expression data sets:  Prostate (Singh et al., 2002) and Neuroblastoma (Oberthuer et al., 2006). The prostate data set contains $p=12601$ genes with $n=136$ samples including 59 positives and 77 negatives. The neuroblastoma  (NB) data set, available from the MicroArray Quality Control phase-II (MAQC-II) project (MAQC Consortium, 2010), consists of gene expression profiles for $p = 10707$ genes from 239 patients (49 positives and 190 negatives) of the German Neuroblastoma Trials NB90-NB2004 with the 3-year event-free survival (3-year EFS) information available. See those references for more detailed description of the data sets.

We fit the logistic regression model with SICA implemented with ICA algorithm (Fan and Lv, 2011). Before applying the regularization method, we exploited the sure independence screening approach to reduce the dimensionality. The random permutation idea (Fan et al., 2011) was applied to determine the threshold for marginal screening. After the screening step, the numbers of retained variables are  430 (prostate) and 2778 (neuroblastoma), respectively. We then chose the final model using those six model selection criteria. Moreover, we randomly split the data into training (80\%) and testing  (20\%) sets for 100 times, and reported the median test classification error rate along with the median model size in Table \ref{tb::real}.

From Table \ref{tb::real}, for the prostate data set the best criterion appears to be $\mbox{GBIC}_p$-L, which has the smallest test classification error rate. For the neuroblastoma data set, if we only look at the median test classification error rate,  $\mbox{GBIC}_p$-L again has the best performance with a small model size. It is worth noting that $\mbox{GBIC}_p$ leads to the most parsimonious model, with median model size 3, at the expense of slightly increasing the test classification error rate. From the results of real examples, it is evident that by taking into account the effect of model misspecification, the performance of the original model selection criteria can be improved in general. This is important since the true model structure is generally unavailable to us in real applications. Our results suggest that the term involving model misspecification in the asymptotic expansions is usually nonnegligible for model selection.

\begin{table}\caption{Results for Prostate and Neuroblastoma data sets.\label{tb::real}}
\begin{center}

\begin{tabular}{c|cccccc}
\hline
\multicolumn{7}{c}{Median classification error rate (in percentage) with RSD in parentheses}\\
\hline
Data set & AIC & BIC & GAIC & GBIC & $\mbox{GBIC}_p$-L & $\mbox{GBIC}_p$ \\
\hline
Prostate &19(9) & 15(6) & 15(9) & 15(9) & 13(9) & 15(10)  \\
NB &18(5) & 18(5) & 18(3) & 18(3) & 18(5) & 19(5)  \\
\hline
\multicolumn{7}{c}{Median model size with RSD in parentheses}\\
\hline
Prostate &15.0(3.7) & 8.5(4.5) & 3.0(1.5) & 6.0(3.7) & 6.0(3.7) & 5.0(3.0)  \\
NB &27.0(3.0) & 26.0(2.2) & 8.5(3.7) & 6.0(3.7) & 5.0(3.0) & 3.0(2.2)  \\
\hline
\end{tabular}

\end{center}

\end{table}

\section{Discussion} \label{sec::discussions}
Despite the rich literature on model selection, the general case of model misspecification in high dimensions is less well studied. Our work has investigated the problem of model selection in high-dimensional misspecified models and characterized the impact of model misspecification. The newly suggested information criterion $\mbox{GBIC}_p$ involving a logarithmic factor of the dimensionality in penalizing model complexity has been shown to perform well in high-dimensional settings. Moreover, we have established the consistency of the covariance contrast matrix estimator that captures the effect of model misspecification in the general setting.

The $\log p$ term in $\mbox{GBIC}_p$ is adaptive to high dimensions. In the setting of correctly specified models, Fan and Tang (2013) showed that such a term is necessary for the model selection consistency of information criteria when the dimensionality diverges fast with the sample size. It would be interesting to study the optimality of those different information criteria under model misspecification. It would also be interesting to investigate model selection principles in more general high-dimensional misspecified models such as the additive models and survival models. These problems are beyond the scope of the current paper and are interesting topics for future research.

\appendix

\section{Proofs of some main results} \label{SecA}
This appendix presents the proofs of Theorems \ref{Thm5} and \ref{Thm3}--\ref{Thm4}. To save space, the proofs of all other theorems and technical lemmas are included in the Supplementary Material. For notational simplicity, throughout the proofs we may specify the orders of different quantities without stating the exact constants, and use the notation $\by$ for observed response and $\bY$ for random response interchangeably when it is convenient.


\subsection{Proof of Theorem \ref{Thm5}} \label{SecA.6}
Throughout this proof $\|\cdot\|_2$ denotes the Euclidean norm of a given vector. The main idea of the proof is to obtain a probabilistic lower bound for the event $\{\hbbeta_n \in N_n(\delta_n)\}.$ To accomplish that we first consider an event that is a subset of this event and calculate the probabilistic lower bound for the smaller event.

Recall that $b''(\theta)$ is continuous and bounded away from $0$, and $\bX$ is of full column rank $d$. Recall the definition $N_n(\delta_n) = \{\bbeta \in \mR^d: \|(n^{-1} \bB_n)^{1/2} (\bbeta - \bbeta_{n, 0})\|_2 \leq (n/d)^{-1/2} \delta_n\}$, and let $\partial N_n(\delta_n)$ denote the boundary of this neighborhood. Since $N_n(\delta_n)$ is compact, $\ell_n(\by, \cdot)$ a continuous strictly concave function, whenever the event
\begin{equation}\label{eq:thm1:proof:event:Q_n}
Q_n = \left\{\ell_n(\by, \bbeta_{n, 0}) > \max_{\bbeta \in \partial N_n(\delta_n)} \ell_n(\by, \bbeta)\right\}
\end{equation}
occurs, $\hbbeta_n$ will be in $N_n(\delta_n)$. The strict concavity of the log-likelihood function follows from the positive definiteness of $\bA_n(\bbeta) = \bX^T \Sig(\bX \bbeta) \bX$, which is the negative of the Hessian of the log-likelihood. This property entails that on the event $Q_n$, the global maximizer $\hbbeta_n$ must belong to the interior of the neighborhood $N_n(\delta_n)$. Hereafter we condition on the event $Q_n$ defined in \eqref{eq:thm1:proof:event:Q_n}. The technical arguments that follow herein, in order to prove that $Q_n$ holds with significant probability, require delicate analyses due to growing dimensionality $d$.


Applying Taylor's expansion to the log-likelihood function $\ell_n(\by, \cdot)$ around $\bbeta_{n, 0}$, we obtain
\[ \ell_n(\by, \bbeta) - \ell_n(\by, \bbeta_{n, 0}) = (\bbeta - \bbeta_{n, 0})\t \bPsi_n(\bbeta_{n, 0}) - \frac{1}{2} (\bbeta - \bbeta_{n, 0})\t \bA_n(\bbeta_*) (\bbeta - \bbeta_{n, 0}), \]
where $\bbeta_*$ is on the line segment joining $\bbeta$ and $\bbeta_{n, 0}$ and $\bPsi_n(\bbeta_{n,0}) = \bX^T[\by - \bmu(\bX\bbeta_{n,0})]$. By letting
$\bu = \delta_n^{-1}  d^{-1/2} \bB_n^{1/2} (\bbeta - \bbeta_{n, 0})$,
the above Taylor's expansion can be rewritten as
\begin{equation}\label{eq:thm1:proof:diff likelihood}
 \ell_n(\by, \bbeta) - \ell_n(\by, \bbeta_{n, 0}) =d^{1/2} \delta_n  \bu\t \bB_n^{-1/2} \bPsi_n(\bbeta_{n, 0}) - d\delta_n^2 \bu\t \bV_n(\bbeta_*) \bu/2,
 \end{equation}
where $\bV_n(\bbeta) = \bB_n^{-1/2}\bA_n(\bbeta)\bB_n^{-1/2}$.

From the definition of $\bu$, $\bbeta \in \partial N_n(\delta_n)$ is equivalent to $\|\bu\|_2 = 1$, and $\bbeta \in \partial N_n(\delta_n)$ implies $\bbeta_* \in N_n(\delta_n)$ since $N_n(\delta_n)$ is convex.
Also it is clear that
\begin{equation}\label{eq:thm1:proof:sup 1}
 \max_{\|\bu\|_2 = 1} \bu\t \bB_n^{-1/2} \bPsi_n(\bbeta_{n, 0}) = \|\bB_n^{-1/2} \bPsi_n(\bbeta_{n, 0})\|_2.
 \end{equation}
From Condition \ref{condition:min_ev_B_n and min_ev_V}, for $n$ sufficiently large, $\min_{\bbeta \in N_n(\delta_n)} \lambda_{\min} \left\{\bV_n(\bbeta)\right\} > c_1n^{-r}$ where $0 < r < 1/4$.
Using this condition and since $\bbeta_* \in N_n(\delta_n)$, it holds that
\begin{equation}\label{eq:thm1:proof:sup 2}
\min_{\|\bu\|_2 = 1} \bu\t \bV_n(\bbeta_*) \bu \geq \min_{\bbeta \in N_n(\delta_n)} \lambda_{\min} \left\{\bV_n(\bbeta)\right\} > c_1n^{-r}.
\end{equation}
Hence by combining \eqref{eq:thm1:proof:sup 1}--\eqref{eq:thm1:proof:sup 2} and taking a supremum on the boundary $\partial N_n(\delta_n)$ in \eqref{eq:thm1:proof:diff likelihood} we derive
\begin{align}\label{eq:thm1:proof:beta_distance}
\max_{\bbeta \in \partial N_n(\delta_n)} \ell_n(\by, \bbeta) - \ell_n(\by, \bbeta_{n, 0}) <& d^{1/2}\delta_n  [\|\bB_n^{-1/2} \bPsi_n(\bbeta_{n, 0})\|_2 \\
&- 2^{-1}c_1 n^{-r}d^{1/2} \delta_n]. \nonumber
\end{align}
By (\ref{eq:normal:equation}), we have $\bX^T[E\by - \bmu(\bX\bbeta_{n,0})]=0$. Hence, $\bPsi_n(\bbeta_{n,0}) = \bX^T[\by - \bmu(\bX\bbeta_{n,0})] = \bX^T(\by-E\by)$.
Denote by $\bW=\bB_n^{-1/2} \bPsi_n(\bbeta_{n, 0})=\bB_n^{-1/2}\bX^T(\by-E\by)$. Notice that $E\bW=\bzero$ and $\cov(\bW)=\bB_n^{-1/2}\cov(\bPsi_n(\bbeta_{n,0})) \bB_n^{-1/2} = \bB_n^{-1/2}\bB_n\bB_n^{-1/2}=I_d$.

Clearly the left hand side of \eqref{eq:thm1:proof:beta_distance} is negative with probability given by
\begin{align}\label{eq:thm1:proof:W-bound}
  P\{\|\bW\|_2 \leq  2^{-1} c_1n^{-r} d^{1/2}\delta_n\}.
\end{align}
From the expression of $\bW$, we have
\begin{align*}
 \|\bW\|_2^2 = & (\by-E\by)\t\bX\bB_n^{-1}\bX\t(\by-E\by)\\
=& [(\by-E\by)\t\cov(\by)^{-1/2}] [ \cov(\by)^{1/2}\bX\bB_n^{-1}\bX\t\cov(\by)^{1/2}]\\
& \cdot [\cov(\by)^{-1/2}(\by-E\by)],
\end{align*}
where $\cdot$ denotes product. Denote by $\bR=\cov(\by)^{1/2}\bX\bB_n^{-1}\bX\t\cov(\by)^{1/2}$ and $\bq=\cov(\by)^{-1/2}(\by\\-E\by)$.
It is easy to check that $\bR^2=\bR$. Therefore, $\bR$ is a projection matrix with rank $\tr(\bR)=d$.
In addition, we have $E\bq=\bzero$ and $\cov(\bq)=I_n$.

We now decompose $\|\bW\|_2^2$ into two terms, the summations of the diagonal entries and the off-diagonal entries, respectively,
\begin{align}
\label{eq:thm1:proof:decomposition}\|\bW\|_2^2=\sum_{i=1}^n r_{ii}q_i^2 + \sum_{1 \leq i\neq j \leq n} r_{ij}q_iq_j,
\end{align}
where $r_{ij}$ denotes the $(i,j)$-entry of $\bR$.
Next we obtain probabilistic bounds for each of the two terms.


From the sub-Gaussian tail condition for $\bq$ in Condition \ref{cond1}, there exists some positive constant $H$ such that for any $t\geq0$,
  \begin{align}
    P(|q_i|> t)\leq H \exp(-t^2/H),
  \end{align}
  for all $1 \leq i \leq n$.

Thus for any $t\geq0$, it holds that
\begin{align}\label{eq:thm1:proof:bound:q_i}
  P\left\{\bigcap_{i=1}^n \{q_i^2\leq t^2\} \right\} &\geq  1-\sum_{i=1}^n P\{q_i^2 > t^2\}\geq  1-nH\exp(-t^2/H).
\end{align}
On the event $\displaystyle\bigcap\nolimits_{i=1}^n \{q_i^2 \leq t^2\}$, we can bound the first term of \eqref{eq:thm1:proof:decomposition} as
\begin{align}\label{eq:thm1:proof:decomposition0}
  \sum_{i=1}^n r_{ii}q_i^2\leq t^2\tr(\bR) =dt^2.
\end{align}

Denote by $\bR_D$ a diagonal matrix with diagonal entries $r_{ii}$. As a result, we observe that $\displaystyle\sum\nolimits_{1 \leq i\neq j \leq n} r_{ij}q_iq_j = \bq\t (\bR - \bR_D) \bq$. It is easy to see that $E \left[\bq\t (\bR - \bR_D) \bq\right] = 0$. We will use a version of the Hanson-Wright inequality (see, e.g., Theorem 1.1 of Rudelson and Vershynin, 2013) to obtain the concentration bound of the quadratic form $\bq\t (\bR - \bR_D) \bq$. But we first start with some notation and preparation.

Let $\| \xi \|_{\psi_2}$ denote the sub-Gaussian norm of a sub-Gaussian random variable $\xi$ defined as
$\| \xi \|_{\psi_2} = \sup_{m \geq 1} \left\{ m^{-1/2}(E | \xi |^m)^{1/m} \right\}$.
From Condition \ref{cond1},  that is, the condition on sub-Gaussian tails, we derive
\begin{align*}
E|q_i|^m&=m\int_{0}^{\infty} x^{m-1} P(|q_i|\geq x) dx\leq Hm\int_{0}^{\infty} x^{m-1} \exp(-x^2/H)dx\\
&=(Hm/2)H^{m/2}\int_{0}^{\infty}u^{m/2-1}\exp(-u)du\\
&=(Hm/2)H^{m/2}\Gamma(m/2)\leq (Hm/2)H^{m/2}(m/2)^{m/2},
\end{align*}
where the last line follows directly from the definition of the Gamma function. Taking the $m$-th root, we have
\[
(E|q_i|^m)^{1/m}  \leq (Hm/2)^{1/m} H^{1/2}(m/2)^{1/2}.
\]
Rewriting after bounding $(1/2)^{(1/m) + (1/2)}$ by 1, we obtain
\[
m^{-1/2} (E|q_i|^m)^{1/m} \leq m^{1/m}H^{(1/2) + (1/m)} \leq e^{1/e} (H^{3/2} \vee 1)
\]
since $m \geq 1$. Therefore, it holds that $\| q_i \|_{\psi_2} \leq c_4$ for all $i$, where $c_4 = e^{1/e} (H^{3/2} \vee 1)$.

We now need bounds on the operator and Frobenius norms of $\bR - \bR_D$. Denote $\| \cdot \|_2$ and $\| \cdot \|_F$ as the matrix operator and Frobenius norms, respectively. Note that $\|\bR\|_2 = 1$ and  $\|\bR\|_F^2 = \tr(\bR^2) = \tr(\bR) = d$. Thus using the fact $\sum_{i\neq j}r_{ij}^2\leq d$, we obtain $\|\bR - \bR_D\|_F^2 \leq d$. Since $|r_{ii}| \leq 1$, we further obtain $\|\bR - \bR_D\|_2 \leq \|\bR\|_2 + \|\bR_D\|_2 \leq 2$. Thereby, a direct application of the Hanson-Wright inequality yields
\begin{align} \label{eq:HW}
P&\left\{\left|\sum_{i\neq j} r_{ij}q_{i}q_{j}\right| > dt^2\right\}
	=P\left\{|\bq^T(\bR - \bR_D)\bq| > dt^2\right\} \\
  	&\leq  2\exp\left\{-c_5 \min \left(\frac{d^2t^4}{c_4^4 \|\bR - \bR_D\|_F^2}, \frac{dt^2}{c_4^2\|\bR - \bR_D\|_2}\right)\right\}\nonumber\\
  	&\leq  2\exp\left\{-c_5 \min \left(\frac{dt^4}{c_4^4}, \frac{dt^2}{2c_4^2}\right)\right\}\nonumber\\
  	&\leq  2 \exp\{-c_6 dt^2\} \nonumber
\end{align}
for any $t > c_4/\sqrt{2}$, where $c_5$ and $c_6$ are some positive constants. To ensure $t > c_4/\sqrt{2}$, we choose $\delta_n = n^{r} (c_0 \log n)^{1/2}$ for some constant $c_0 > 8c_1^{-2}H$ and $t=2^{-3/2}c_1n^{-r}\delta_n$. Therefore, the probability bound \eqref{eq:HW} holds for large enough $n$.

Combining \eqref{eq:thm1:proof:bound:q_i} and \eqref{eq:HW}, with probability at least
$1-nH\exp(-t^2/H)\\-2\exp(-c_6dt^2)$, we have
\begin{align*}
  \|\bW\|_2^2&\leq  \sum_i r_{ii}q_i^2+\left|\sum\nolimits_{i\neq j} r_{ij}q_iq_j\right|  2 dt^2.
\end{align*}
In view of our choice of $t$, it holds that $(2dt^2)^{1/2} = 2^{-1} c_1n^{-r} d^{1/2}\delta_n$ and thus
\begin{align*}
P\{\|\bW\|_2  \leq 2^{-1} c_1n^{-r} d^{1/2}\delta_n\} &\geq 1-nH\exp(-t^2/H)-2\exp(-c_6dt^2)\\
& \geq1 - O(n^{-\alpha}),
\end{align*}
where $\alpha = [(c_1^2c_0/(8H)-1)] \wedge (c_1^2c_6c_0/8)$. Note that $\alpha > 0$ since $c_0 > 8c_1^{-2}H$. This leads to
 \begin{equation}\label{eq:thm1:proof:bound:probnew}
 P(Q_n) \geq 1-O(n^{-\alpha}).
 \end{equation}
The positive constant $\alpha$ can be large if $c_0$ in $\delta_n$ is chosen to be large. From Condition \ref{condition:min_ev_B_n and min_ev_V}, $\lammin(\bB_n) \rightarrow \infty$ at a faster rate than $d\delta_n^2$. Then we have the consistency $\hbbeta_n - \bbeta_{n,0} = o_P(1)$.


\subsection{Proof of Theorem \ref{Thm3}} \label{SecA.4}
Define $\mathcal{E}=\{\widehat\bbeta_n\in N_n(\delta_n)\}$, where  $\widehat\bbeta_n$ stands for the QMLE. Note that $\mathcal{E}$ does depend on $n$, but for simplicity of notation we will omit the subscript $n$ in sequel. To establish this theorem we require a possibly dimension dependent bound on the quantity $\| n^{-1/2} \bX \hbbeta_n \|_2$. The need for bounding the specified quantity, particularly with growing dimensionality, can be intuitively understood by trying to put some restriction on the parameter space. This is analogous to the case of penalized likelihood.

Recall the neighborhood $M_n(\alpha_1) = \{\bbeta \in \mR^d: \|\bX \bbeta \|_{\infty} \leq \alpha_1 \log n\}$, where $\alpha_1$ is some positive constant satisfying $\alpha_1 < \alpha/2 - 1$. One way of bounding the quantity $\| n^{-1/2} \bX \hbbeta_n \|_2$ is to restrict the QMLE $\hbbeta_n$ on the set $M_n(\alpha_1$).  As mentioned in Theorem \ref{Thm5}, the constant $\alpha$ can be large if $c_0$ is chosen to be large, which ensures that $\alpha_1$ is positive. From Condition \ref{condition:region}, $N_n(\delta_n) \subset M_n (\alpha_1)$ for all sufficiently large $n$ to ensure that conditional on $\mathcal{E}$, the restricted MLE coincides with its unrestricted version. However, this condition is very mild in the sense that the constant $\alpha_1$ can be chosen as large as desired to make $M_n(\alpha_1)$ large enough, whereas the neighborhood $N_n(\delta_n)$ is asymptotically shrinking. Hereafter in this proof $\hbbeta_n$ will be referred to as the restricted MLE, unless specified otherwise.

Recall that $\eta_n(\bbeta)=E\ell_n(\widetilde \by,\bbeta)$, where $\widetilde \by$ is an independent copy of $\by$. In the GLM setup, we have $\ell_n(\widetilde \by,\bbeta) = \widetilde \by^T\bX\bbeta - \bone^T \bb(\bX\bbeta)$ and $\eta_n(\bbeta)=(E\widetilde \by^T)\bX\bbeta-\bone^T\bb (\bX\bbeta)$.

\smallskip

\textbf{Part 1: Expansion of $E\eta_n(\widehat\bbeta_n)$. }
We approach the proof by splitting $E\eta_n(\widehat\bbeta_n)$ in the region $\mathcal{E}$ and its complement, that is,
\begin{align} \label{split}
E\eta_n(\widehat\bbeta_n)
&=E\{\eta_n(\widehat\bbeta_n)1_{\mathcal{E}}\}+E\{\eta_n(\widehat\bbeta_n)1_{\mathcal{E}^c}\}\\ \nonumber
&= E\{\eta_n(\widehat\bbeta_n)1_{\mathcal{E}}\} +  E\{[(E\tilde \by)^T(\bX\widehat\bbeta_n)-\bone^T\bb(\bX\widehat\bbeta_n)]1_{\mathcal{E}^c}\},
\end{align} where the second equality follows from the definition of $\eta_n(\cdot)$.

We aim to show that the second term on the right hand side of \eqref{split} is $o(1)$. Performing componentwise Taylor's expansion of $\bb(\cdot)$ around $\bzero$ and evaluating at $\bX\hbbeta_n$, we obtain
$\bb(\bX\hbbeta_n) = \bb(\bzero) + b'(0)\bX\hbbeta_n + \br$,
 where $\br = (r_1, \cdots, r_n)^T$ with $r_i =  2^{-1} b''((\bX\bbeta_i^*)_i) (\bX\hbbeta_n)_i^2$ and $\bbeta_1^*, \cdots ,\bbeta_n^*$ lying in the line segment joining $\hbbeta_n$ and $\bzero$. Recall that $\hbbeta_n $ is the constrained MLE here, $EY_i^2$ is bounded uniformly in $i$ and $n$, and $b''(\cdot) = O(n^{\alpha_1})$ uniformly in its argument. The condition on $b''(\cdot)$ can be much weakened in many cases including linear and logistic regression models. This condition also accommodates Poisson regression where $b''(\theta) = \exp(\theta)$ for $\theta \in \mathbb{R}$ since $b(\theta) = \exp(\theta)$.
Then it follows that
\begin{align}\label{orderonethm5}
E\{|(E\widetilde \by)^T\bX\widehat\bbeta_n-\bone^T\bb(\bX\widehat\bbeta_n)|1_{\mathcal{E}^c}\} 
& \leq O\{n \log n +  n + n^{1 + \alpha_1} (\log n)^2\}P(\mathcal{E}^c) \\
&\leq O\{n^{2(\alpha_1 + 1)}\}P(\mathcal{E}^c) = o(1) \nonumber
\end{align}
for sufficiently large $n$. The last inequality follows from the fact that $\alpha > 2(\alpha_1 +1)$ and we recall that $P(\mathcal{E}^c) = O(n^{-\alpha})$. To verify the orders, we note that the four bounds $|(E\widetilde \by)^T\bX\widehat\bbeta_n| \leq n \max_{1 \leq i \leq n} (E y_i^2)^{1/2} \alpha_1 \log n $, $|\bone^T \bb(\bzero)| = O(n)$, $|b'(0) \bone^T \bX \hbbeta_n| \leq O(1) n \alpha_1 \log n$, and $|\bone^T \br| \leq n \max_{1 \leq i \leq n} |r_i| \leq $ \\ $n O(n^{\alpha_1})(\alpha_1 \log n)^2$.

On the event $\mathcal{E}$, we first expand $\eta_n(\bbeta)$ around $\bbeta_{n,0}$. By the definition of  $\bbeta_{n,0}$, $\eta_n(\bbeta)$ attains its maximum at $\bbeta_{n,0}$. By Taylor's expansion of $\eta_n(\cdot)$ around $\bbeta_{n,0}$ and evaluating at $\hbbeta_n$, we derive
\begin{align} \label{basicexpansion}
\eta_n(\hbbeta_n) &= \eta_n(\bbeta_{n,0}) - \frac{1}{2} (\hbbeta_n - \bbeta_{n,0})^T \bA_n(\bbeta^*) (\hbbeta_n - \bbeta_{n,0}) \\
&= \eta_n(\bbeta_{n,0}) - \frac{1}{2} (\hbbeta_n - \bbeta_{n,0})^T \bA_n (\hbbeta_n - \bbeta_{n,0}) - \frac{s_n}{2}\nonumber \\
&= \eta_n(\bbeta_{n,0}) - \frac{1}{2} \bv_n^T [(\bC_n^{-1})^T \bA_n \bC_n^{-1}] \bv_n - \frac{s_n}{2},  \nonumber
\end{align}
where $\bA_n (\cdot) = - \partial^2 \ell_n(\by, \cdot)/\partial \bbeta^2$, $\bA_n = \bA_n(\bbeta_{n,0})$, $s_n = (\hbbeta_n - \bbeta_{n,0})^T [\bA_n(\bbeta^*) -   \bA_n] (\hbbeta_n - \bbeta_{n,0})$, $\bv_n = \bC_n (\hbbeta_n - \bbeta_{n,0})$, and $\bbeta^*$ is on the line segment joining $\bbeta_{n,0}$ and $\hbbeta_n$. Then it follows that
\begin{align} \label{orderapprox}
\left|s_n 1_\mathcal{E}\right| &= \left|(\hbbeta_n - \bbeta_{n,0})^T (\bA_n(\bbeta^*) -   \bA_n) (\hbbeta_n - \bbeta_{n,0})\right| 1_\mathcal{E} \\
&= \left |[\bB_n^{1/2} (\hbbeta_n - \bbeta_{n,0})]^T [\bV_n(\bbeta^*) - \bV_n] [\bB_n^{1/2} (\hbbeta_n - \bbeta_{n,0})]\right| 1_\mathcal{E} \nonumber \\
& \leq \  \|\bV_n(\bbeta^*) - \bV_n \|_2 \delta_n^2 d 1_\mathcal{E}, \nonumber
\end{align}
where $\bV_n(\cdot) = \bB^{-1/2}\bA_n(\cdot)\bB_n^{-1/2}$ and $\bV_n = \bV(\bbeta_{n,0})$. Note that on the event $\mathcal{E}$, by the convexity of the neighborhood $N_n(\delta_n)$ we have $\bbeta^* \in N_n(\delta_n)$. From Condition \ref{condition:widetilde_V_n}, $\max_{\bbeta_1, \cdots, \bbeta_d \in N_n(\delta_n)} \|\widetilde{\bV}_n(\bbeta_1, \cdots \bbeta_d) - \bV_n\|_2 = O(d^{1/2}n^{-1/2})$. Therefore we deduce that $E(s_n 1_{\mathcal{E}})$ is of order $O(d^{3/2}n^{-1/2}\delta_n^2) =  o(1)$, which follows from (\ref{order_AN}) in the proof of Theorem \ref{Thm6}.

From (\ref{decom_AN}) in the proof of Theorem \ref{Thm6}, we have the decomposition $\bv_n = \bu_n + \bw_n$ with $\bu_n = \bB_n^{-1/2} \bX\t (\by - E \by)$ and
\[ \bw_n = -\left[\widetilde{\bV}_n(\bbeta_1, \cdots, \bbeta_d) - \bV_n\right] \left[\bB_n^{1/2} (\hbbeta_n - \bbeta_{n, 0})\right]. \] For simplicity of notation, denote by $\bR_n = (\bC_n^{-1})^T \bA_n \bC_n^{-1}$. Recall that $\bC_n = \bB_n^{-1/2}\bA_n$. With some calculations we obtain
\begin{align*}
E(\bu_n^T  \bR_n \bu_n) &= E\{(\by-E\by)\t\bX\bA_n^{-1} \bX\t(\by-E\by)\} \\ &= E\{\tr(\bA_n^{-1}\bX\t(\by-E\by)(\by-E\by)\t\bX)\} = \tr(\bA_n ^{-1} \bB_n).
\end{align*}
Note that $E (\bu_n^T  \bR_n \bu_n 1_{\mathcal{E}}) = E(\bu_n^T  \bR_n \bu_n) - E (\bu_n^T  \bR_n \bu_n 1_{\mathcal{E}^c})$.
From Theorem \ref{Thm5}, we have $P(\mathcal{E}^c) \rightarrow 0$ as $n \rightarrow \infty$. Let $\mu_n = \tr(\bA_n^{-1} \bB_n) \vee 1$ ensuring that this quantity is bounded away from zero. We will apply Vitali's convergence theorem to show that $E (\bu_n^T  \bR_n \bu_n 1_{\mathcal{E}^c}) = o(\mu_n)$. To establish uniform integrability we use the following lemma, the proof of which has been provided in Appendix \ref{SecB} in Supplementary Material.
\begin{lemma}\label{UI}
For some constant $\gamma>0$, $\sup_n E |(\bu_n^T \bR_n \bu_n)/\mu_n |^{1+\gamma} < \infty$.
\end{lemma}
This leads to $E (\bu_n^T  \bR_n \bu_n 1_{\mathcal{E}^c}) = o(\mu_n)$. Hence we have
\[
\frac{1}{2} E (\bu_n^T  \bR_n \bu_n 1_{\mathcal{E}}) = \frac{1}{2}\tr(\bA_n ^{-1} \bB_n) + o(\mu_n). \nonumber
\]

It remains to show that
\begin{equation} \label{eq:thm 5 part 1 last}
E[(\bw_n^T\bR_n\bw_n + 2\bw_n^T\bR_n\bu_n)1_{\mathcal{E}}] = o(\mu_n).
\end{equation} Note that on the event $\mathcal{E}$, we have
\begin{align}
\bw_n^T\bR_n\bw_n &= \| \bR_n^{1/2} \bw_n \|_2^2 \nonumber \leq \| \widetilde{\bV}_n - \bV_n \|_2^2 \delta_n^2d \tr(\bA_n^{-1} \bB_n). \nonumber
\end{align}
In view of the assumption $\max_{\bbeta_1, \cdots, \bbeta_d \in
  N_n(\delta_n)} \|\widetilde{\bV}_n(\bbeta_1, \cdots, \bbeta_d) -
\bV_n\|_2 = O(d^{1/2}n^{-1/2})$, it holds that $E(\bw_n^T\bR_n\bw_n1_{\mathcal{E}}) = o(\mu_n)$. For the cross term $\bw_n^T\bR_n\bu_n$, applying the Cauchy-Schwarz inequality yields
\begin{align}
|E(\bw_n^T\bR_n\bu_n 1_{\mathcal{E}})| &\leq E(\|\bR_n^{1/2} \bw_n\|_2^21_{\mathcal{E}})^{1/2} E(\| \bu_n^T\bR_n^{1/2} \|_2^2)^{1/2} \nonumber \\
&\leq E[\| \widetilde{\bV}_n - \bV_n \|_2 1_{\mathcal{E}} \delta_n d^{1/2} \tr(\bA_n^{-1} \bB_n)], \nonumber
\end{align}
which entails that $E(\bw_n^T\bR_n\bu_n1_{\mathcal{E}}) = o(\mu_n)$. Note that $E\{|\eta_n(\bbeta_{n,0})| 1_{\mathcal{E}^c}\}$ is of order $o(1)$ by similar calculations as in (\ref{orderonethm5}). Thus combining (\ref{split}) -- (\ref{eq:thm 5 part 1 last}) yields
$E\{\eta_n(\hbbeta_n)\}= \eta_n(\bbeta_{n,0}) - \frac{1}{2} \tr(\bA_n^{-1}\bB_n) + o(\mu_n)$.

\smallskip

\textbf{Part 2: Expansion of $E\ell_n(\by,\bbeta_{n,0})$.} Similarly we expand $\ell_n(\by, \cdot)$ around $\hbbeta_n$ and evaluate at $\bbeta_{n,0}$. From Condition \ref{condition:region}, $N_n(\delta_n) \subset M_n (\alpha_1)$ for sufficiently large $n$, we see that $\bbeta_{n,0} \in M_n (\alpha_1)$. On the event $\mathcal{E}$, since $\ell_n(\by,\cdot)$ attains its maximum at the restricted MLE $\hbbeta_n$, we have
\begin{align}
\ell_n(\by, \bbeta_{n,0}) &= \ell_n(\by, \hbbeta_n) - \frac{1}{2} (\hbbeta_n - \bbeta_{n,0})^T \bA_n(\bbeta^*)(\hbbeta_n - \bbeta_{n,0}) \\
&= \ell_n(\by, \hbbeta_n) - \frac{1}{2} (\hbbeta_n - \bbeta_{n,0})^T \bA_n (\hbbeta_n - \bbeta_{n,0}) -\frac{s_n}{2} \nonumber \\
&= \ell_n(\by, \hbbeta_n) - \frac{1}{2} \bv_n^T [(\bC_n^{-1})^T \bA_n \bC_n^{-1}] \bv_n \nonumber - \frac{s_n}{2}.\label{ell_expansion}
\end{align}
Then similarly as in Part 1, we can obtain
\[E\{\ell_n(\by, \bbeta_{n,0})1_{\mathcal{E}}\}  = E\{\ell_n(\by,\hbbeta_n)1_{\mathcal{E}}\} - \frac{1}{2} \tr(\bA_n^{-1}\bB_n) + o(\mu_n).\] If we can show that $E\{|\ell_n(\by, \bbeta_{n,0})|1_{\mathcal{E}^c}\}$ and $E\{|\ell_n(\by,\hbbeta_n)|1_{\mathcal{E}^c}\}$ are both of order $o(1)$, then we obtain the desired asymptotic expansion
\[
E\{\eta_n(\hbbeta_n)\} = E\{\ell_n(\by,\hbbeta_n)\} - \tr(\bA_n^{-1}\bB_n) + o(\mu_n).
\]

To see why $E\{|\ell_n(\by, \bbeta_{n,0})|1_{\mathcal{E}^c}\}$ is of order $o(1)$, we derive
\begin{align*}
E&\{|\by^T\bX\bbeta_{n,0} - \bone^T \bb(\bX\bbeta_{n,0})|1_{\mathcal{E}^c}\} \\
&\leq O(n \log n) P(\mathcal{E}^c)^{1/2}  +  O\{n + n \log n+ n^{2 + \alpha_1} (\log n)^2\}P(\mathcal{E}^c) \\
&\leq O\{n^{2(\alpha_1 + 1)}\}P(\mathcal{E}^c) = o(1),
\end{align*}
similarly as in (\ref{orderonethm5}) and using $E[|\by^T\bX\bbeta_{n,0}| 1_{\mathcal{E}^c}] \leq E[|\by^T\bX\bbeta_{n,0}|^2]^{1/2} P(\mathcal{E}^c)^{1/2}$. Similarly we can also show that $E\{|\ell_n(\by, \hbbeta_n)|1_{\mathcal{E}^c}\}$ is of order o(1). The only difference in the above derivation is to bound $\|\bX \hbbeta_n\|_{\infty}$ instead of $\|\bX \bbeta_{n,0}\|_{\infty}$, which holds from the definition of the restricted QMLE. This concludes the proof.

\subsection{Proof of Theorem \ref{Thm4}} \label{SecA.5} In view of the expansions of GAIC, GBIC, and $\mbox{GBIC}_p$, we need to show that $\log |\hH_n| = \log |\bH_n| + o_P(1)$ and $\tr(\hH_n) = \tr(\bH_n) + o_P(1)$. To establish this we show that $\hH_n = \bH_n + o_P(1/d)$, where the $o_P(\cdot)$ denotes the convergence in probability of the matrix operator norm.

Let $\bM$ be a $d \times d$ square matrix. Denote by $\overline{\tr} (\bM) = \tr (\bM)/d$ the normalized trace and $\rho(\bM) = \max_{1 \leq k \leq d} \{ |\lambda_k(\bM)| \}$ the spectral radius. Then we have
\begin{align*}
|\tr(\hH_n) - \tr(\bH_n)| &= d |\overline{\tr} (\hH_n - \bH_n)| \\
&\leq d \rho(\hH_n - \bH_n) = d \| \hH_n - \bH_n \|_2 = o_P(1),
\end{align*}
where $\| \cdot \|_2$ denotes the matrix operator norm. The equality of the spectral radius and the operator norm follows from the symmetry of the matrix $\hH_n - \bH_n$. Similarly define the normalized log determinant, that is, $\overline{\log}\,|\bM| = (\log |\bM|)/d$ for any arbitrary matrix $\bM$. Denote $\lambda_k(\cdot)$ as the eigenvalues arranged in the increasing order. Then we have
\begin{align} \label{log_det_order}
|\log|\hH_n| - \log |\bH_n| | & \leq d | \overline \log \, |\hH_n| - \overline \log \, |\bH_n| | \nonumber \\ & \leq d \max_{1 \leq k \leq d} |\log \lambda_k(\hH_n) - \log \lambda_k(\bH_n) | \nonumber \\ & \leq d \max_{1 \leq k \leq d} \log \left\{1 + \left|\frac{\lambda_k(\hH_n)}{\lambda_k(\bH_n)} - 1\right|\right\}.
\end{align}
Recall that we assume that the smallest and largest eigenvalues of both $n^{-1}\bB_n$ and $n^{-1}\bA_n$ are bounded away from 0 and $\infty$. It then follows that $\lambda_k(\bH_n) =O(1)$ and $\lambda_k^{-1}(\bH_n) =O(1)$ uniformly for all $1 \leq k \leq d$. An application of Weyl's theorem shows that
\[
|\lambda_k(\hH_n) - \lambda_k(\bH_n)| \leq \rho(\hH_n - \bH_n)
\]
for each $k$. We have $\rho(\hH_n - \bH_n) = \| \hH_n - \bH_n \|_2 = o_P(1/d)$. Hence the right hand side of (\ref{log_det_order}) is $o_P(1)$.

Now we proceed to show that $\hH_n = \bH_n + o_P(1/d)$. It suffices to prove that $n^{-1}\hA_n = n^{-1} \bA_n + o_P(1/d)$ and $n^{-1}\hB_n = n^{-1} \bB_n + o_P(1/d)$. 
We use the following properties of the operator norm (Horn and Johnson, 1985): $\| (I_d-\bM)^{-1} \|_2 \leq 1/(1-\| \bM \|_2)$ if $\| \bM \|_2 < 1$, $\| \bM \bN \|_2 \leq \| \bM \|_2 \| \bN \|_2$, and $\| \bM + \bN \|_2 \leq \| \bM \|_2 + \| \bN \|_2$, where $\bM$ and $\bN$ are $d \times d$ matrices. 
To see the sufficiency note that
\begin{align*}
&(n^{-1} \hA_n)^{-1} (n^{-1} d \hB_n) - (n^{-1} \bA_n)^{-1} (n^{-1} d \bB_n) \\=& (n^{-1} \hA_n)^{-1} (n^{-1} d \hB_n)- (n^{-1} \hA_n)^{-1} (n^{-1} d \bB_n) + (n^{-1} \hA_n)^{-1} (n^{-1} d \bB_n) \\&- (n^{-1} \bA_n)^{-1} (n^{-1} d \bB_n).
\end{align*}

Then the desired result $\hH_n = \bH_n + o_P(1/d)$ can be obtained by repeated application of the above properties of the operator norm.

\smallskip

\textbf{Part 1:  Prove $n^{-1}\hA_n = n^{-1} \bA_n + o_P(1/d)$}.
From Theorem \ref{Thm5} we have, $\|(n^{-1}\bB_n)^{1/2} (\hbbeta_n - \bbeta_{n,0})\|_2 = O_P\{(n/d)^{-1/2}\delta_n\}$, which along with the assumption that the smallest eigenvalue of $n^{-1}\bB_n$ is bounded away from 0 entails $\hbbeta_n = \bbeta_{n,0} + O_P\{(n/d)^{-1/2}\delta_n\}$.
Then it follows from the Lipschitz assumption for $n^{-1} \bA_n(\bbeta)$ from Condition \ref{condition:lipschitz} in the neighborhood $N_n(\delta_n)$ and Theorem \ref{Thm5} that $n^{-1}\hA_n = n^{-1} \bA_n + o_P(1/d)$, which holds for our choice of $d = o\{n^{(1-4r)/3} (\log n)^{-2/3}\}$ and $\delta_n$.

\smallskip

\textbf{Part 2: Prove $n^{-1}\hB_n = n^{-1} \bB_n + o_P(1/d)$}. We first split $n^{-1} \hB_n$ as
\[
n^{-1} \hB_n = n^{-1} \bX\t \diag\left\{\left[\by - \bmu(\bX
  \hbbeta_n)\right] \circ \left[\by - \bmu(\bX
  \hbbeta_n)\right]\right\} \bX = \bG_1 + \bG_2 + \bG_3,
\]
where
\begin{align*}
\bG_1 & = n^{-1} \bX\t \diag\{(\by - \bmu(\bX \bbeta_{n, 0})) \circ (\by - \bmu(\bX \bbeta_{n, 0}))\} \bX, \\
\bG_2 & = 2 n^{-1} \bX\t \diag\{(\by - \bmu(\bX \bbeta_{n, 0})) \circ [\bmu(\bX \bbeta_{n, 0}) - \bmu(\bX \hbbeta_n)]\} \bX, \\
\bG_3 & = n^{-1} \bX\t \diag\{[\bmu(\bX \hbbeta_n) - \bmu(\bX \bbeta_{n, 0})] \circ [\bmu(\bX \hbbeta_n) - \bmu(\bX \bbeta_{n, 0})]\} \bX.
\end{align*}
We will state two lemmas before proceeding with the proof. Define the sub-exponential norm of a sub-exponential random variable $\xi$ as
\[
\| \xi \|_{\psi_1} = \sup_{m \geq 1} \left\{ m^{-1}(E | \xi |^m)^{1/m} \right\}.
\]

\begin{lemma}\label{lemma::7} For independent sub-Gaussian random variables $\{y_i\}_{i=1}^n$, we have that $q_i^2 = (y_i - Ey_i)^2/\var(y_i)$ is sub-exponential with norm bounded by $2c_4^2$, where $c_4$ is as defined in the proof of Theorem \ref{Thm5}. Moreover, the following Bernstein-type tail probability bound holds
\begin{equation} \label{subexp}
P\left\{|\Sigma_{i=1}^n a_iq_i^2 - E[\Sigma_{i=1}^n a_iq_i^2]| \geq t \right\} \leq 2 \exp \left[ - c_{10} \min \left( \frac{t^2}{4c_4^4 \|\ba\|_2^2}, \frac{t}{2c_4^2 \|\ba\|_{\infty}} \right) \right] \nonumber
\end{equation}
for $\ba \in \mathbb{R}^n$, $t \geq 0$, and $c_{10} > 0$.
\end{lemma}


\begin{lemma}\label{lemma::8} For independent sub-Gaussian random variables $\{y_i\}_{i=1}^n$ with $q_i = \{\var(y_i)\}^{-1/2}\\(y_i - Ey_i)$, the following tail probability bound holds
\begin{equation} \label{subexp}
P\left\{|\Sigma_{i=1}^n a_iq_i| \geq t \right\} \leq e \exp \left(-\frac{c_{11}t^2}{c_4^2 \|\ba\|_2^2} \right) \nonumber
\end{equation}
for $\ba \in \mathbb{R}^n$, $t \geq 0$, and $c_{11} > 0$.
\end{lemma}
Lemma \ref{lemma::7} follows from Lemma 5.14 and Proposition 5.16 of Vershynin (2012). Note that here we define the sub-exponential random variable as the square of a sub-Gaussian random variable and the bound on the norm follows by our previous observation that $\|q_i\|_{\Psi_2} \leq c_4$ in the proof of Theorem \ref{Thm5}. Lemma \ref{lemma::8} rephrases Proposition 5.10 of Vershynin (2012) for the case where $\|q_i\|_{\Psi_2} \leq c_4$.

Further split $\bG_1$ as $\bG_1 = \bG_{11} + \bG_{12} + \bG_{13}$ where
\begin{align*}
\bG_{11} & = n^{-1} \bX\t \diag\{(\by - E\by) \circ (\by - E\by)\} \bX, \\
\bG_{12} & = 2 n^{-1} \bX\t \diag\{(\by - E\by) \circ [E\by - \bmu(\bX \bbeta_{n, 0})]\} \bX, \\
\bG_{13} & = n^{-1} \bX\t \diag\{[E\by - \bmu(\bX \bbeta_{n, 0})] \circ [E\by - \bmu(\bX \bbeta_{n, 0})]\} \bX.
\end{align*}
Note that $E\bG_{11} = n^{-1}\bB_n$ and $\bG_{11} = n^{-1} \Sigma_{i=1}^n \{ \bx_i \bx_i^T [y_i - Ey_i]^2\} = \Sigma_{i=1}^n \bA_i q_i^2$, where $\bA_i = n^{-1} \var (y_i) \bx_i \bx_i^T$. Then it holds that for any positive $t$,
\begin{align}\label{eq::G11-bound}
P(\|\bG_{11} - E\bG_{11}\|_2 \geq t) & \leq P(\|\bG_{11} - E\bG_{11}\|_F \geq t)  \nonumber\\
&\leq d^2 \max_{1 \leq j,k \leq d} P(|\bG_{11}^{jk} - E\bG_{11}^{jk}| \geq t/d),
\end{align}
where $\| \cdot \|_F$ denotes the matrix Frobenius norm and $\bG_{11}^{jk}$ denotes the $(j,k)$ entry of $\bG_{11}$. Recall from Condition  \ref{condition:lipschitz} that $\| \bX \|_{\infty} = O(n^{\alpha_2})$ with $0 \leq \alpha_2 < r$. Define $a_i^{jk} = n^{-1} \var(y_i) x_{ij} x_{ik}$ and $\ba^{jk} = (a_1^{jk}, \cdots ,a_n^{jk})^T$. We have $\|\ba^{jk}\|_2^2 = O(n^{-1}n^{4 \alpha_2})$.
Then combining \eqref{eq::G11-bound} with Lemma \ref{lemma::7}, we deduce
\begin{align*}
P(d\|\bG_{11} - E\bG_{11}\|_2 \geq t) &\leq d^2 \max_{1 \leq j,k \leq d} P(|\bG_{11}^{jk} - E\bG_{11}^{jk}| \geq t/d^2) \\& \leq 2d^2 \exp \{ - c_{12}t^2n^{1-4\alpha_2}/d^4 \}
\end{align*}
for some constant $c_{12} > 0$. Note that $d = o\{n^{(1-4r)/4}\}$, we obtain $\bG_{11} = E\bG_{11} + o_P(1/d)$.

By Condition \ref{condition:bias_variance} and Lemma \ref{lemma::8}, we have 
\begin{align*}
P(d\|\bG_{12}\|_2 \geq t) &\leq d^2P(|\bG_{12}^{jk}| \geq t/d^2)  \leq ed^2 \exp \{ - c_{13}t^2n^{1 - 4\alpha_2 + (1 - \alpha_3)/2}/d^4 \},
\end{align*}
where $c_{13} > 0$ is some constant. Hence from  $d = o\{n^{(1-4r)/4}\}$ and $0 \leq \alpha_3 \leq 4(r - \alpha_2)$, we have $\bG_{12} = o_P(1/d)$.

To show that $\bG_{13} = o(1/d)$, we derive
\begin{align*}
\|\bG_{13}\|_2^2 &\leq \|n^{-1} \Sigma_{i=1}^n \{ \bx_i \bx_i^T [Ey_i - [\bmu(\bX\bbeta_{n,0})]_i]^2\}\|_F^2 \\ & = \Sigma_{1 \leq j,k \leq d} [\Sigma_{i=1}^n a_{i}^{jk} [Ey_i - [\bmu(\bX\bbeta_{n,0})]_i]^2/\var(y_i)]^2 \\ & \leq  \Sigma_{i=1}^n \{[Ey_i - [\bmu(\bX\bbeta_{n,0})]_i]^2/\var(y_i)\}^2 \Sigma_{1 \leq j,k \leq d} \|\ba^{jk}\|_2^2,
\end{align*}
where the last step follows from the component-wise Cauchy-Schwarz inequality. From Condition \ref{condition:bias_variance},  $\bG_{13} = o(1/d)$. Combining the above derivations yields
$\bG_{1} = E\bG_{1} + o_P(1/d) = n^{-1}\bB_n + o_P(1/d)$.
To see that $\bG_2 = o_P(1/d)$, note that $(\by - \bmu(\bX \bbeta_{n, 0}))_i = (y_i - Ey_i) + (Ey_i  - [\bmu(\bX \bbeta_{n, 0})]_i)$ and apply similar arguments as above. By the Lipschitz Condition \ref{condition:lipschitz} in the neighborhood $N_n(\delta_n)$, we have $\bG_3 = o_P(1/d)$, which completes the proof.



\newpage
\setcounter{page}{1}
\setcounter{section}{0}
\setcounter{equation}{0}

\renewcommand{\theequation}{A.\arabic{equation}}
\setcounter{equation}{0}
\appendix

\begin{center}{\bf \Large Supplementary Material to ``Model Selection in High-Dimensional Misspecified Models"}

\bigskip

Pallavi Basu, Yang Feng and Jinchi Lv
\end{center}

\noindent This Supplementary Material contains the proofs of Theorems \ref{Thm6} and \ref{Thm1}--\ref{Thm2}, and technical lemmas, as well as additional tables from Section \ref{sec::simulation}.

\smallskip

\appendix
\setcounter{page}{1}
\setcounter{section}{1}
\renewcommand{\theequation}{A.\arabic{equation}}
\setcounter{equation}{0}

\section{Proofs of Additional Theorems}\label{SecA1}

\subsection{Proof of Theorem \ref{Thm6}} \label{SecA.7}
Recall that $\bC_n = \bB_n^{-1/2} \bA_n$. To establish the asymptotic normality of the QMLE $\hbbeta_n$, we prove the following
  \begin{align}
    \bD_n \bC_n(\widehat\bbeta_n-\bbeta_{n,0}) \toD N(\bzero,I_m),
  \end{align}
  for any $m\times d$ matrix $\bD_n$ such that $\bD_n\bD_n^T =I_m$ with $m$ fixed.
From the score equation we have $\bPsi(\hbbeta_n) = \bX^T[\by - \bmu(\bX\hbbeta_n)] = 0$. From (\ref{eq:normal:equation}), it holds that $\bX^T[E\by - \bmu(\bX\bbeta_{n,0})] = 0$.  For any $\bbeta_1,
\cdots, \bbeta_d \in \mathbb{R}^d$, denote by $\widetilde{\bA}_n(\bbeta_1, \cdots,
\bbeta_d)$ a $d \times d$ matrix with $j$-th row the corresponding row of $\bA_n(\bbeta_j)$ for each $j = 1, \cdots, d$, and matrix-valued function $\widetilde{\bV}_n(\bbeta_1, \cdots \bbeta_d) = \bB_n^{-1/2} \widetilde{\bA}_n(\bbeta_1, \cdots, \bbeta_d) \bB_n^{-1/2}$. Assuming the differentiability of $\bPsi(\cdot)$ and applying the mean-value theorem componentwise around $\bbeta_{n,0}$, we obtain
\begin{align*}
\bzero & = \bPsi_n(\hbbeta_n) = \bPsi_n(\bbeta_{n, 0}) - \widetilde{\bA}_n(\bbeta_1, \cdots, \bbeta_d) (\hbbeta_n - \bbeta_{n, 0}) \\
& = \bX\t (\by - E \by) - \widetilde{\bA}_n(\bbeta_1, \cdots, \bbeta_d) (\hbbeta_n - \bbeta_{n, 0}),
\end{align*}
where each of $\bbeta_1, \cdots, \bbeta_d$ lies on the line segment joining $\hbbeta_n$ and $\bbeta_{n, 0}$. It follows from this expansion that
\begin{equation} \label{decom_AN}
\bC_n (\hbbeta_n - \bbeta_{n, 0}) = \bu_n + \bw_n,
\end{equation}
where $\bu_n = \bB_n^{-1/2} \bX\t (\by - E \by)$ and
\[ \bw_n = -\left[\widetilde{\bV}_n(\bbeta_1, \cdots, \bbeta_d) - \bV_n\right] \left[\bB_n^{1/2} (\hbbeta_n - \bbeta_{n, 0})\right], \]
where $\bV_n =\bV_n(\bbeta_{n, 0}) = \bB_n^{-1/2} \bA_n \bB_n^{-1/2}$. Therefore we have
\[ \bD_n\bC_n (\hbbeta_n - \bbeta_{n, 0}) = \bD_n\bu_n + \bD_n\bw_n.\]
By the Cram\'{e}r-Wold theorem, it suffices to show that for any unit vector $\ba\in \mR^m$, $\ba^T\bD_n\bC_n (\hbbeta_n - \bbeta_{n, 0})\toD N(0,1)$. Further by Slutsky's lemma, it is sufficient to show that
$\ba^T\bD_n\bu_n\toD N(0,1)$ and $\ba^T\bD_n\bw_n =o_P(1)$ for any unit vector $\ba$.

\smallskip

\textbf{Part 1 (Asymptotic normality of $\ba^T\bD_n\bu_n$): }
We will build on the conditions required to apply the Lyapunov central limit theorem (CLT). For an arbitrary unit vector $\ba \in \mR^m$, consider the asymptotic distribution of
\[ v_n = \ba\t \bD_n\bu_n = \ba\t \bD_n\bB_n^{-1/2} \bX\t (\by - E \by) = \sum_{i = 1}^n z_i, \]
where $z_i = \ba\t \bD_n\bB_n^{-1/2} \bx_i (y_i - E y_i)$, $i = 1, \cdots, n$, and $\bX = (\bx_1, \cdots, \bx_n)^T$. Since $z_i$'s are independent and have mean zero, we derive
\begin{align*}
\var(v_n) &= \sum_{i = 1}^n \var(z_i) = \ba\t \bD_n\bB_n^{-1/2} \bX^T\cov(\by)\bX \bB_n^{-1/2} \bD_n\t\ba \\
&= \ba\t \bD_n\bB_n^{-1/2} \bB_n \bB_n^{-1/2} \bD_n\t\ba = 1.
\end{align*}
From Condition \ref{condition:Yi_deviation and sub-gaussian}, we have $\max_{1\leq i\leq n} E |y_i - E y_i|^3 \leq M$ for some positive constant $M$ and $\sum_{i = 1}^n (\bx_i\t \bB_n^{-1}\bx_i)^{3/2} = o(1)$. Then an application of the Cauchy-Schwarz inequality yields
\begin{align*}
\sum_{i = 1}^n E |z_i|^3 & = \sum_{i = 1}^n |\ba\t \bD_n\bB_n^{-1/2} \bx_i|^3 E |y_i - E y_i|^3 \leq M \sum_{i=1}^n |\ba\t \bD_n\bB_n^{-1/2} \bx_i |^3 \\
& \leq M \sum_{i=1}^n \|\bD_n\t\ba\|_2^3\|\bB_n^{-1/2} \bx_i\|_2^3 = M \sum_{i=1}^n (\bx_i\t \bB_n^{-1} \bx_i)^{3/2}\rightarrow 0,
\end{align*}
noting that $\|\bD_n\t\ba\|_2^2 = \ba^T\bD_n\bD_n^T\ba = \ba^T I_m \ba = 1$. Therefore by applying Lyapunov's CLT, we obtain
\[ \ba\t \bD_n\bu_n = \sum_{i = 1}^n z_i \toD N(0, 1). \]

\textbf{Part 2 (To show $\ba^T\bD_n\bw_n$ is o(1) in probability): }
Conditional on the event $\{\hbbeta_n \in N_n(\delta_n)\}$ and using the fact that $\|\bD_n\t\ba\|_2 = 1$, we have
\begin{align*}
  |\ba^T\bD_n\bw_n|&\leq  \|\bD_n^T \ba\|_2\|\bw_n\|_2 \leq  \|\bw_n\|_2\\
  &\leq \|\widetilde\bV_n-\bV_n\|_2\|\bB_n^{1/2}(\widehat\bbeta_n-\bbeta_{n,0})\|_2\\
  &\leq  \|\widetilde\bV_n-\bV_n\|_2 d^{1/2}\delta_n,
\end{align*}
where the last step follows from the definition of the neighborhood $N_n(\delta_n) = \{\bbeta \in \mR^d: \|(n^{-1} \bB_n)^{1/2} (\bbeta - \bbeta_{n, 0})\|_2 \leq (n/d)^{-1/2} \delta_n\}$ and given that $\{\hbbeta_n \in N_n(\delta_n)\}$.
From Condition \ref{condition:widetilde_V_n}, $\max_{\bbeta_1, \cdots, \bbeta_d \in
  N_n(\delta_n)} \|\widetilde{\bV}_n(\bbeta_1, \cdots, \bbeta_d) -
\bV_n\|_2 \\= O(d^{1/2}n^{-1/2}) \leq O(dn^{-1/2}\delta_n)$.
Again conditional on the event $\{\hbbeta_n \in N_n(\delta_n)\}$ and noticing that each $\bbeta_j$ defined previously for $1 \leq j \leq d$ lies in $N_n(\delta_n)$ due to its convexity, it holds that
\begin{align} \label{order_AN}
  |\ba\t\bD_n\bw_n|
  =  O(d^{3/2}n^{-1/2}\delta_n^2) =  o(1),
\end{align}
where we choose $\delta_n = n^r (c_0\log n)^{1/2}$ as in the proof of Theorem \ref{Thm5} and $d = o\{n^{(1-4r)/3}\\ (\log n)^{-2/3}\}$ with $0 \leq r < 1/4$. Since the event $\{\hbbeta_n \in N_n(\delta_n)\}$ holds with probability tending to 1, $\ba\t\bD_n\bw_n = o_P(1)$. Also note that the convergence to zero in probability is uniform in $\ba$ and $\bD_n$. Therefore, combining parts 1 and 2 finishes the proof.

\subsection{Proof of Theorem \ref{Thm1}} \label{SecA.2}
Throughout the proof we condition on the event $\widetilde Q_n =\{\hbbeta_n\in N_n(\delta_n)\}$, where $N_n(\delta_n)=\{\bbeta\in \mR^d: \|(n^{-1}\bB_n)^{1/2}(\bbeta-\bbeta_{n,0})\|_2\leq (n/d)^{-1/2} \delta_n\}$, $\bB_n = \bX\t \cov(\bY) \bX$, and $\hbbeta_n$ is the unrestricted MLE. From Theorem \ref{Thm5} we have shown that as $n \rightarrow \infty$, $$P(\widetilde Q_n)\rightarrow 1.$$

Recall from (\ref{eq:ellnstar}) that $\ell_n^*(\by,\bbeta) = \ell_n(\by, \bbeta) - \ell_n(\by, \hbbeta_n)$. Then the maximum value zero of this function is attained at $\bbeta=\widehat\bbeta_n$. It follows from (\ref{090}) that
\[ \partial^2 \ell^*_n(\by, \bbeta)/\partial \bbeta^2 = -\bA_n(\bbeta),\] where $\bA_n(\bbeta) = \bX^T \bSigma(\bX\bbeta)\bX$. By Taylor's expansion of the likelihood function $\ell_n(\by, \cdot)$ around $\widehat\bbeta_n$ in the new neighborhood $\widetilde N_n(\delta_n)=\{\bbeta\in \mR^d: \|(n^{-1}\bB_n)^{1/2}(\bbeta-\widehat\bbeta_n)\|_2\leq (n/d)^{-1/2} \delta_n\}$, we derive
\begin{align} \label{026}
\ell^*_n(\by, \bbeta) &= \frac{1}{2} (\bbeta-\widehat\bbeta_n)\t \left[\partial^2 \ell^*_n(\by, \bbeta_*)/\partial \bbeta^2\right] (\bbeta-\widehat\bbeta_n) \\  &= -\frac{n}{2}\bdelta\t\bV_n(\bbeta_*)\bdelta, \nonumber
\end{align}
where $\bbeta_*$ lies on the line segment joining $\bbeta$ and $\widehat\bbeta_n$,
$\bdelta=n^{-1/2}\bB_n^{1/2}(\bbeta-\widehat\bbeta_n)$, and $\bV_n(\bbeta)=\bB_n^{-1/2}\bA_n(\bbeta)\bB_n^{-1/2}$. Since $\widehat\bbeta_n \in \widetilde N_n(\delta_n)$, by the convexity of the neighborhood $\widetilde N_n(\delta_n)$ we have  $\bbeta_* \in \widetilde N_n(\delta_n)$. Also note that conditional on the event $\widetilde Q_n$, it holds that $\widetilde N_n(\delta_n) \subset N_n(2\delta_n)$. We define $$\rho_n(\delta_n) = \max_{\bbeta \in N_n(2\delta_n)} \max \{ | \lambda_{\min}(\bV_n(\bbeta) - \bV_n) |, | \lambda_{\max} (\bV_n(\bbeta) - \bV_n) | \}$$ with $\bV_n = \bV_n(\bbeta_{n,0})$. Using Taylor's expansion (\ref{026}) over the region $\widetilde N_n(\delta_n)$, we obtain
\begin{equation} \label{027}
q_1 (\bbeta)1_{\widetilde N_n(\delta_n)} (\bbeta) \leq -n^{-1} \ell^*_n(\by, \bbeta) 1_{\widetilde N_n(\delta_n)} (\bbeta) \leq q_2 (\bbeta) 1_{\widetilde N_n(\delta_n)} (\bbeta),
\end{equation}
where $q_1 (\bbeta) = \frac{1}{2} \bdelta^T[\bV_n - \rho_n(\delta_n)I_d]\bdelta$ and $q_2 (\bbeta) = \frac{1}{2} \bdelta^T[\bV_n + \rho_n(\delta_n)I_d]\bdelta$.

Define $U_n(\bbeta) = \exp\left[n^{-1} \ell^*_n(\by, \bbeta)\right]$ which takes values in the interval $[0,1]$ by definition. From Condition \ref{condition:min_ev_B_n and min_ev_V}, for $n$ large, $\min_{\bbeta \in N_n(\delta_n)} \lambda_{\min} \left\{\bV_n(\bbeta)\right\} > c_1n^{-r}$ with $0 < r < 1/4$ and $\rho_n(\delta_n) = o\{n^{-(1-r)/3}\}$. Since $\bbeta_{n,0}$ belongs to $N_n(\delta_n)$, this assumption yields $\rho_n(\delta_n) \leq \lambda_{\min}(\bV_n)/2$ for sufficiently large $n$. To see this, note that since $(1-r)/3 > r$ we have $\rho_n(\delta_n) n^r = o(1)$ whereas $\lambda_{\min}(\bV_n) n^r > c_1$. Consider the linear transformation $h(\bbeta) = (n^{-1}\bB_n)^{1/2}\bbeta$. For sufficiently large $n$, we obtain
\begin{align} \label{028}
E_{\mu_{\mathfrak{M}}} [e^{-nq_2(\bbeta)}  1_{\widetilde N_n(\delta_n)}(\bbeta)] &\leq E_{\mu_{\mathfrak{M}}} [U_n(\bbeta)^n 1_{\widetilde N_n(\delta_n)}(\bbeta)] \\ &\leq E_{\mu_{\mathfrak{M}}} [e^{-nq_1(\bbeta)}  1_{\widetilde N_n(\delta_n)}(\bbeta)], \nonumber
\end{align}
where $\mu_{\mathfrak{M}}$ denotes the prior distribution on $h(\bbeta) \in \mathbb{R}^d$ for model $\mathfrak{M}$. Before proceeding with the proof we state a few useful lemmas. The proofs of these lemmas are elaborated in Appendix \ref{SecA}. From Condition \ref{condition:prior density and rho_n}, the prior density relative to the Lebesgue measure $\mu_0$ on $\mathbb{R}^d$, $\pi(h(\bbeta)) = \frac{d\mu_{\mathfrak{M}}}{d\mu_0} (h(\bbeta))$,  satisfies
\begin{equation}\label{condonpi}
\inf_{\bbeta \in N_n(2\delta_n)} \pi (h(\bbeta)) \geq c_2 \text{  and  } \sup_{\bbeta \in \mathbb{R}^d} \pi (h(\bbeta)) \leq c_3,
\end{equation}
where $c_2$ and $c_3$ are some positive constants.

\begin{lemma}\label{lemma::2}
Under (\ref{condonpi}), for $j = 1,2,$ we have
\begin{equation}\label{029}
c_2 \int_{\bdelta \in \mathbb{R}^d} e^{-nq_j} 1_{\widetilde N_n (\delta_n)} d\mu_0 \leq E_{\mu_{\mathfrak{M}}} \left[e^{-nq_j} 1_{\widetilde N_n (\delta_n)}\right] \leq c_3 \int_{\bdelta \in \mathbb{R}^d} e^{-nq_j} 1_{\widetilde N_n (\delta_n)} d\mu_0.
\end{equation}
\end{lemma}

\begin{lemma}\label{lemma::3}
Conditional on the event $\widetilde Q_n$, for sufficiently large $n$ we have
\begin{align} \label{030}
E_{\mu_{\mathfrak{M}}} [U_n(\bbeta)^n 1_{\widetilde N_n^c (\delta_n)}] &\leq \exp \{- [\kappa_n - \rho_n(\delta_n)/2]  d \delta_n^2\} \\
&\leq \exp[-(\kappa_n/2) d \delta_n^2], \nonumber
\end{align}
where $\kappa_n = \lambda_{\min}(\bV_n)/2$.
\end{lemma}

\begin{lemma}\label{lemma::4} It holds that
\begin{equation}
\int_{\bdelta \in \mathbb{R}^d} e^{-n q_1} d \mu_0 = \left( \frac{2 \pi}{n} \right)^{d/2} |\bV_n - \rho_n(\delta_n) I_d|^{-1/2}
\end{equation}
and
\begin{equation}
\int_{\bdelta \in \mathbb{R}^d} e^{-n q_2} d \mu_0 = \left( \frac{2 \pi}{n} \right)^{d/2} |\bV_n + \rho_n(\delta_n) I_d|^{-1/2}.
\end{equation}
\end{lemma}

\begin{lemma}\label{lemma::5} For $j = 1,2$, it holds that
\begin{equation}\label{035}
\int_{\bdelta \in \mathbb{R}^d} e^{-nq_j} 1_{\widetilde N_n^c (\delta_n)} d\mu_0 \leq \left( \frac{2\pi}{n\kappa_n} \right)^{d/2} \exp \left[-(\sqrt{\kappa_nd\delta_n^2} - \sqrt{d})^2/2\right]
\end{equation}

\end{lemma}

Now we proceed with the proof. From Condition \ref{condition:min_ev_B_n and min_ev_V}, $\delta_n = n^{r} (c_0 \log n)^{1/2}$. Then the expression in (\ref{030}) converges to zero faster than any polynomial rate in $n$. Let us rewrite the right hand side of (\ref{035}) as $$\exp \left\{ -\frac{d}{2} (\sqrt{\kappa_n\delta_n^2} - 1)^2 + \frac{d}{2} [\log(2 \pi) - \log(n \kappa_n)] \right\},$$ which converges to zero faster than any polynomial rate in $n$. From Condition \ref{condition:min_ev_B_n and min_ev_V}, $d = o\{n^{(1-4r)/3} (\log n)^{-2/3}\}$ with $0 < r < 1/4$.

Then it follows that
\begin{align*}
|\bV_n \pm \rho_n(\delta_n)I_d|^{-1/2} &= |\bV_n|^{-1/2} |\bI_d \pm \rho_n(\delta_n)\bV_n^{-1}|^{-1/2} \\
&= |\bV_n|^{-1/2} \{ 1 + O[\rho_n(\delta_n)\tr(\bV_n^{-1})]\} \\
&= |\bV_n|^{-1/2} \{1 + O[\rho_n(\delta_n)d\lambda_{\min}^{-1}(\bV_n)]\} \\
&= |\bV_n|^{-1/2}[1 + o(1)].
\end{align*}
Combining Lemmas \ref{lemma::2}--\ref{lemma::5} yields
\begin{align*}
\log E_{\mu_{\mathfrak{M}}} [U_n(\bbeta)^n] &= \log \left\{ \left( \frac{2\pi}{n} \right)^{d/2} |\bV_n|^{-1/2} [1 + o(1)] \right\} + \log c_n \\
&= - \frac{\log n}{2} d + \frac{1}{2} \log |\bA_n^{-1} \bB_n| + \frac{\log(2\pi)}{2}d + \log c_n + o(1),
\end{align*}
where $c_n \in [c_2,c_3]$. This completes the proof.

\subsection{Proof of Theorem \ref{Thm2}} \label{SecA.3}
The proof follows from the proof of Theorem \ref{Thm1} and Part 1, that is, expansion of $E \eta_n(\hbbeta_n)$, of the proof of Theorem \ref{Thm3}.

\section{Proofs of Lemmas} \label{SecB}
Lemmas \ref{lemma::7} and \ref{lemma::8}  have been discussed in the paragraph following them. The proofs of Lemmas \ref{lemma::2}--\ref{lemma::4} can be found in Lv and Liu (2014).

\subsection{Proof of Lemma \ref{UI}} From the expression of $\bu_n^T  \bR_n \bu_n$, we have
\begin{align*}
\bu_n^T  \bR_n \bu_n = & (\by-E\by)\t\bX\bA_n^{-1}\bX\t(\by-E\by)\\
=& [(\by-E\by)\t\cov(\by)^{-1/2}] [ \cov(\by)^{1/2}\bX\bA_n^{-1}\bX\t\cov(\by)^{1/2}]\\
& \cdot [\cov(\by)^{-1/2}(\by-E\by)].
\end{align*}
Denote by $\bS_n=\cov(\by)^{1/2}\bX\bA_n^{-1}\bX\t\cov(\by)^{1/2}$ and $\bq=\cov(\by)^{-1/2}(\by-E\by)$. We decompose $\bu_n^T  \bR_n \bu_n$ into two terms, the summations of the diagonal entries and the off-diagonal entries, respectively,
\begin{align}
\bu_n^T  \bR_n \bu_n=\sum_{i=1}^n s_{ii}q_i^2 \ + \sum_{1 \leq i\neq j \leq n} s_{ij}q_iq_j, \nonumber
\end{align}
where $s_{ij}$ denotes the $(i,j)$-entry of $\bS_n$. Then we have
\begin{align*}
E(\bu_n^T  \bR_n \bu_n)^2 \ =& \ \sum_{i=1}^n s_{ii}^2E(q_i^4) \ + \ \sum_{1 \leq i\neq j \leq n} s_{ii} s_{jj} E(q_i^2)E(q_j^2) \ \\ &+ \ \sum_{1 \leq i\neq j \leq n} s_{ij}^2E(q_i^2)E(q_j^2). \nonumber
\end{align*}
Using the sub-Gaussian norm bound $c_4$, both quantities $E(q_i^4)$ and $E(q_i^2)E(q_j^2)$ can be uniformly bounded by a common constant. Hence \[ E(\bu_n^T  \bR_n \bu_n)^2 \leq O(1) \cdot \{[\tr(\bS_n)]^2  + \tr(\bS_n^2)\}. \] Since $\bS_n$ is positive semidefinite it holds that $\tr(\bS_n^2) \leq [\tr(\bS_n)]^2$. Finally noting that $\tr(\bS_n) = \tr(\bA_n^{-1}\bB_n)$, we see that $\sup_n E |(\bu_n^T \bR_n \bu_n)/\mu_n |^{1+\gamma} < \infty$ for $\gamma = 1$.

%

\subsection{Proof of Lemma \ref{lemma::5}} From the definition of $q_j(\bbeta)$ for $j = 1, 2$, we derive
\begin{align}
\exp(-nq_j) =& \exp(- (n/2) \bdelta^T[\bV_n \pm \rho_n(\delta_n)I_d] \bdelta) \nonumber \\
\leq& \exp(- n (\kappa_n - \rho_n(\delta_n)/2) \bdelta^T \bdelta) \nonumber \\
\leq& \exp(- (n \kappa_n)/2 \bdelta^T \bdelta).
\end{align}
Then we have
\begin{align*}
\int_{\bdelta \in \mathbb{R}^d} e^{-nq_j} 1_{\widetilde N_n^c (\delta_n)} d\mu_0 \leq& \int_{\bdelta \in \mathbb{R}^d} e^{- \frac{n \kappa_n}{2} \bdelta^T \bdelta} 1_{\widetilde N_n^c (\delta_n)} d\mu_0 \\
=& \left( \frac{2\pi}{n\kappa_n} \right)^{d/2} P (\| (n\kappa_n)^{-1/2} \bZ \|_2^2 \geq (n/d)^{-1} \delta_n^2) \\
=& \left( \frac{2\pi}{n\kappa_n} \right)^{d/2} P (\|\bZ \|_2^2 \geq \kappa_n d \delta_n^2),
\end{align*}
where $\bZ \sim N(\bzero, I_d)$.

Using the chi-square tail bound, that is, for any positive $x$ it is known that $P (\|\bZ \|_2^2 - d \geq 2\sqrt{dx} + 2x) \leq \exp(-x)$ and after minor modification it holds that $P (\|\bZ \|_2^2 \geq (\sqrt{d} + \sqrt{2x})^2) \leq \exp(-x)$. With this observation, define $x = (\sqrt{\kappa_nd\delta_n^2} - \sqrt{d})^2/2$ and the proof concludes.

\section{Additional Tables}
In Tables \ref{tb::simu-linear-inter-weak-FP-FN}--\ref{tb::simu-logi-inter-FP-FN}, we report additional variable selection results for the three simulation examples in Section \ref{sec::simulation}.

\begin{table}[!ht]
\caption{Example \ref{sec::simu-linear-interaction-weak}. Median false positives  with median false negatives (strong/weak effects) in parentheses when the model is misspecified.
\label{tb::simu-linear-inter-weak-FP-FN}}
\begin{center}

\begin{tabular}{c|cccccc}
\hline
& AIC & BIC & GAIC & GBIC & $\mbox{GBIC}_p$-L & $\mbox{GBIC}_p$ \\
\hline
200 & 24(0/4) & 2(0/5) & 2(0/5) & 1(0/5) & 0(0/5) & 0(0/5) \\
400 & 23(0/4) & 19(0/5) & 3(0/5) & 4(0/5) & 0(0/5) & 0(0/5) \\
1600 & 24(0/5) & 25(0/5) & 4(0/5) & 23(0/5) & 1(0/5) & 0(0/5) \\
3200 & 19(0/5) & 25(0/5) & 4(0/5) & 25(0/5) & 2(0/5) & 0(0/5) \\
\hline
\end{tabular}

\end{center}

\end{table}

\begin{table}[!ht]
\caption{Example \ref{sec::simu-linear-interaction-weak}. Median false positives  with median false negatives (strong/weak effects) in parentheses when the model is correctly specified.
\label{tb::simu-linear-inter-weak-FP-FN-correct}}
\begin{center}

\begin{tabular}{c|cccccc}
\hline
& AIC & BIC & GAIC & GBIC & $\mbox{GBIC}_p$-L & $\mbox{GBIC}_p$ \\
\hline
200 & 73(0/0) & 0(0/4) & 0(0/4) & 0(0/4) & 0(0/4) & 0(0/4) \\
400 & 77(0/0) & 0(0/4) & 0(0/4) & 0(0/4) & 0(0/4) & 0(0/5) \\
1600 & 4(0/2) & 0(0/5) & 0(0/5) & 0(0/5) & 0(0/5) & 0(0/5) \\
3200 & 0(0/5) & 0(0/5) & 0(0/5) & 0(0/5) & 0(0/5) & 0(0/5) \\
\hline
\end{tabular}

\end{center}

\end{table}

\begin{table}[!ht]\caption{Example \ref{sec::simu-multiple-index}. Median false positives with median false negatives in parentheses.
\label{tb::simu-multiple-index-FP-FN}}
\begin{center}
\begin{tabular}{c|cccccc}
\hline
 $p$ & AIC & BIC & GAIC & GBIC & $\mbox{GBIC}_p$-L & $\mbox{GBIC}_p$ \\
\hline
200 & 4(0) & 4(0) & 4(0) & 3(0) & 0(0) & 0(0) \\
400& 5(0) & 5(0) & 5(0) & 5(0) & 4(0) & 0(0) \\
1600& 8(0) & 8(0) & 8(0) & 8(0) & 7(0) & 0(0) \\
3200& 8(0) & 8(0) & 8(0) & 8(0) & 8(0) & 0(0) \\
\hline\end{tabular}

\end{center}

\end{table}

\begin{table}[!ht]\caption{Example \ref{sec::simu-logi-inter}. Median false positives  with median false negatives in parentheses.
\label{tb::simu-logi-inter-FP-FN}}
\begin{center}
\begin{tabular}{c|cccccc}
\hline
 $p$ & AIC & BIC & GAIC & GBIC & $\mbox{GBIC}_p$-L & $\mbox{GBIC}_p$ \\
\hline
200& 29(0) & 1(0) & 11(0) & 1(0) & 0(0) & 0(0) \\
400& 25(0) & 5(0) & 14(0) & 1(0) & 0(0) & 0(0) \\
1600 &19(0) & 18(0) & 14(0) & 3(0) & 1(0) & 0(0) \\
3200& 18(0) & 17(0) & 13(0) & 9(0) & 1(0) & 0(0) \\
\hline
\end{tabular}
\end{center}

\end{table}

\end{document}